\documentclass[reqno,11pt]{amsart}
\usepackage{epsf}\usepackage{amssymb} \usepackage{cite}

\numberwithin{equation}{section}
\newcommand{\beq}{\begin{equation}}
\newcommand{\ee}{\end{equation}}
\newcommand{\bea}{\begin{eqnarray}}
\newcommand{\eea}{\end{eqnarray}}

\usepackage{hyperref}
\def\stackreb#1#2{\ \mathrel{\mathop{#1}\limits_{#2}}}

\newcommand{\CC}{\mathbb C}
\newcommand{\R}{\mathbb R}
\newcommand{\Z}{\mathbb Z}

\newcommand{\T}{\mathbb T}

\begin{document}

\title[Complex and rational hypergeometric functions]
{Complex and rational hypergeometric
\\ functions on root systems}

\author{ G. A. Sarkissian and  V. P. Spiridonov}

\makeatletter
\renewcommand{\@makefnmark}{}
\makeatother

\address{
G. A. Sarkissian: Laboratory of Theoretical Physics, JINR, Dubna, 141980,  Russia; Yerevan Physics Institute,
Alikhanian Br. 2, 0036, Yerevan, Armenia; e-mail: sarkissn@theor.jinr.ru
\vspace{-2mm}}

\address{
V. P. Spiridonov: Laboratory of Theoretical Physics, JINR, Dubna, 141980,  Russia;
National Research University Higher School of Economics, Moscow, Russia; email: spiridon@theor.jinr.ru}

\begin{abstract}
We consider some new limits for the elliptic hypergeometric integrals on root systems.
After the degeneration of elliptic beta integrals of type I and type II for root systems
$A_n$ and $C_n$ to the hyperbolic hypergeometric integrals, we apply the limit
$\omega_1\to - \omega_2$ for their quasiperiods (corresponding to $b\to \textup{i}$
in the two-dimensional conformal field theory) and obtain complex beta integrals in
the Mellin--Barnes representation admitting exact evaluation.
Considering type I elliptic hypergeometric integrals of a higher order
obeying nontrivial symmetry transformations, we derive their descendants to the level
of complex hypergeometric functions and prove the Derkachov--Manashov conjectures
for functions emerging in the theory of non-compact spin chains.
We describe also symmetry transformations for a type II complex hypergeometric function
on the $C_n$-root system related to the recently derived generalized complex Selberg integral.
For some hyperbolic beta integrals we consider a special limit
$\omega_1\to \omega_2$ (or $b\to 1$) and obtain new hypergeometric
identities for sums of integrals of rational functions.
\end{abstract}

\maketitle


\tableofcontents

\section{Introduction}

Ordinary hypergeometric functions define classical special functions which have found very many applications
in theoretical and mathematical physics as well as in the pure mathematical problems \cite{aar}.
The theory of $q$-deformed hypergeometric functions  was developing since the time of Euler, but its
applications in physics emerged only about forty years ago. It is through the mathematical physics
problems that relatively recently the third type of special functions of hypergeometric type,
depending on two basic parameters $q$ and $p$, has been discovered, see e.g. \cite{spi:essays} for a survey.
The most general form of these functions is described by the elliptic hypergeometric integrals \cite{spi:umn}
comprising all previously known ordinary and $q$-hypergeometric functions as special limiting degenerations.
These integrals represent the functions which are transcendental
over the field of elliptic functions. They define meromorphic functions of parameters whose nonsingular
small $p$ expansions yield infinite series with $q$-hypergeometric function coefficients.
By now, the elliptic hypergeometric integrals have found profound applications
in the theory of integrable many-body systems \cite{spi:cs}, the Yang--Baxter equation  \cite{DS:umn}
and, most interestingly, in the quantum field theory \cite{DO,SV}.

Already in the  1960s, an extension of the Euler--Gauss $_2F_1$ hypergeometric function to the field
of complex numbers has been considered by Gelfand, Graev and Vilenkin \cite{GGV}.
Recently this function and the corresponding complex hypergeometric equation were rigorously
investigated from the functional analysis point of view in \cite{MN}.
Although this complex hypergeometric function and its higher order analogues together with
their multivariable extensions have found prominent applications
in the two-dimensional conformal field theory \cite{DF} and in the theory of noncompact spin
chains \cite{DMV2017}, until recently it was not clear where these functions lie in the general
hierarchy of special functions of hypergeometric type. This question was clarified in the
series of works \cite{DMV2018,DSS,DS2017,GSVS1,GSVS3}, where it was rigorously shown that
the complex hypergeometric functions emerge in the Mellin--Barnes
type representation from the hyperbolic hypergeometric integrals which, in turn, appear in a
specific degeneration limit from the elliptic hypergeometric integrals \cite{rai:limits}
(formally such a limit was considered earlier in \cite{vDS3}).
Altogether this gives a unified picture of special functions of the hypergeometric type.

In the present paper we would like to extend limiting considerations of \cite{GSVS1,GSVS3}
from the univariate case to the multiple complex hypergeometric integrals on root systems.
For that we consider elliptic hypergeometric integrals of type I on root systems $A_n$ and $C_n$,
both the beta integrals admitting exact evaluations and higher order integrals obeying nontrivial
symmetry transformations
\cite{spi:essays}. Then we pass to the hyperbolic level
described by infinite contour integrals of combinations of the Faddeev modular quantum dilogarithm
\cite{Fad94,Fad95} called also the hyperbolic gamma function \cite{ruij}.
At the next step of degenerations, we go down to infinite
bilateral sums of the Mellin--Barnes type integrals representing complex hypergeometric
functions in the spirit of Naimark's consideration of $3j$-symbols for the group $SL(2,\mathbb{C})$
\cite{Naimark} and Ismagilov's results on the $6j$-symbols of the same group \cite{Ismag2}
reexamined recently in \cite{DS2017}. In this way we prove the Derkachov--Manashov conjectures
\cite{DM2019} on complex analogues of the symmetry transformations of certain multiple
Gustafson integrals \cite{gus:some1,gus:some2}.

An elliptic analogue of the exactly computable Selberg integral was introduced in \cite{vDS1,vDS2}.
As its special limiting form, a general complex hypergeometric analogue of the Selberg integral
in the Mellin--Barnes representation was recently constructed in \cite{GSVS3}.
Here we extend this result and derive a complex hypergeometric analogue of the symmetry
transformation for a type II elliptic
hypergeometric integral on the root system $C_n$ established by Rains in \cite{rai:trans}.
In the last section we present exact relations for rational functions of hypergeometric
type which arise from a different type of degeneration limit for the hyperbolic hypergeometric integrals
and extend the results of \cite{GSVS2} to the multivariable setting.

Let us remind briefly some structural elements of the theory of complex hypergeometric functions.
Euler's beta integral evaluation formula
\begin{equation}
\int_0^1x^{\alpha-1} (1-x)^{\beta-1}dx =\frac{\Gamma(\alpha)\Gamma(\beta)}{\Gamma(\alpha+\beta)},
\quad \textrm{Re}(\alpha), \textrm{Re}(\beta)>0,
\label{EulerBeta}\end{equation}
where $\Gamma(x)$ is the Euler gamma function,
is an important exact identity which has found numerous applications in mathematics and theoretical
physics, including computations of Feynman integrals. Therefore its generalizations represent
great research interest.
A complex analogue of formula \eqref{EulerBeta} has been  found by Gelfand, Graev and Vilenkin in \cite{GGV}
\begin{equation}\label{CB}
\int_{\mathbb{C}}[w-z_1]^{\alpha-1} [z_2-w]^{\beta-1} \frac{d^2w}{\pi}
=\frac{\Gamma(\alpha)\Gamma(\beta)}{\Gamma(\alpha+\beta)}
\frac{\Gamma(1-\alpha'-\beta')}{\Gamma(1-\alpha')\Gamma(1-\beta')}[z_2-z_1]^{\alpha+\beta-1},
\end{equation}
where $\alpha, \alpha'\in\mathbb{C}$ are complex numbers
such that $\alpha-\alpha' \in\mathbb{Z}$ and the following notation is used
$$
[z]^\alpha:= z^\alpha \bar z^{\alpha'}=(-1)^{\alpha-\alpha'}[-z]^\alpha,\quad
\int_{\mathbb{C}} d^2z:=\int_{\mathbb{R}^2}d(\text{Re}\, z)\, d(\text{Im}\, z),
$$
with $\bar z$ being the complex conjugate of $z$. A drawback of this notation is that
for integer $\alpha$ it looks awkward, e.g.  $[z]^1=z\bar z^{1-N}$ for $\alpha=1,\, \alpha-\alpha'=N$,
i.e. $1'=1-N$. In \eqref{CB} and below it is assumed that in the power exponents one has $1'=1$.

Formula \eqref{CB} becomes maximally close to the Euler's beta integral evaluation in terms of the gamma
function over the field of complex numbers
\begin{equation}
{\bf \Gamma}(\alpha|\alpha'):=\frac{\Gamma(\alpha)}{\Gamma(1-\alpha')}.
\label{Cgamma}\end{equation}
We shall use also a slightly different notation, namely,
$$
{\bf \Gamma}(x,n):={\bf \Gamma}\big( \tfrac{n+\textup{i}x}{2}
| \tfrac{-n+\textup{i} x}{2} \big)
=\frac{\Gamma(\frac{n+\textup{i}x}{2})}{\Gamma(1+\frac{n-\textup{i}x}{2})},
$$
where $x\in \CC$ and $n\in\Z$.
Reflection equations  have the form
\beq
{\bf \Gamma}(\alpha|\alpha') =(-1)^{\alpha-\alpha'}{\bf \Gamma}(\alpha'|\alpha), \qquad
{\bf \Gamma}(x,-n)=(-1)^n{\bf \Gamma}(x,n),
\label{reflCgamma0}\ee
and
\beq
{\bf \Gamma}(\alpha|\alpha'){\bf \Gamma}(1-\alpha|1-\alpha')  =(-1)^{\alpha-\alpha'}, \qquad
{\bf \Gamma}(x,n){\bf \Gamma}(-x-2\textup{i},n)=1.
\label{reflCgamma}\ee
Functional equations
\begin{eqnarray*} &&
{\bf \Gamma}(\alpha+1|\alpha') ={\bf \Gamma}(x-\textup{i},n+1)=\alpha{\bf \Gamma}(\alpha|\alpha'),\quad
\\ &&
{\bf \Gamma}(\alpha|\alpha'+1) ={\bf \Gamma}(x-\textup{i},n-1)=-\alpha' {\bf \Gamma}(\alpha|\alpha')
\end{eqnarray*}
are easily deduced from the definitions $\alpha=\tfrac{n+\textup{i}x}{2},\, \alpha'=\tfrac{-n+\textup{i}x}{2}$
and the standard relation $\Gamma(x+1)=x\Gamma(x)$.

The Stirling formula for the gamma function yields
$$
\frac{ \Gamma(a+z)}{\Gamma(b+z)}=z^{a-b}(1+O(\tfrac{1}{|z|})), \quad |z| \to\infty, \quad |\arg(z)-\pi|>\delta>0.
$$
Consider the following product of two complex gamma functions:
$$
f(z):={\bf \Gamma}(x+R,n+N){\bf\Gamma}(y-R,m-N)=\frac{\Gamma(a+z)\Gamma(c-z)}
{\Gamma(b+{\bar z})\Gamma(d-{\bar z})},
$$
where $R\in \mathbb{R}$,
$$
a=\frac{n+\textup{i}x}{2}, \; b=1+\frac{n-\textup{i}x}{2}, \;
c=\frac{m+\textup{i}y}{2}, \; d=1+\frac{m-\textup{i}y}{2},\; z=\frac{N+\textup{i}R}{2}.
$$
We would like to obtain the large modulus $|z|$ asymptotics of $f(z)$.
Suppose that $0<\arg(z)<\pi$, then for $|z|\to\infty$ we can write
\begin{eqnarray}\nonumber &&
f(z)=\frac{\Gamma(a+z)\Gamma(1-d+{\bar z})}{\Gamma(1-c+z)\Gamma(b+{\bar z})}\frac{\sin[\pi(d-{\bar z})]}{\sin[\pi(c-z)]}
\\  && \makebox[2em]{}
=|z|^{\textup{i}(x+y)-2}e^{\textup{i}[\arg(z)(n+m)+\pi (N- m)]}
(1+O(\tfrac{1}{|z|})).
\label{asympgam}\end{eqnarray}
The asymptotics in the region $-\pi<\arg(z)<0$ is obtained simply by the replacement $\arg(z)\to \arg(z)-\pi$.
In the case $\arg(z)=0$ (i.e. $R=0,\, N>0$), we obtain for large $N$
\[
f\big(\tfrac{N}{2}\big)=\frac{\Gamma(a+\tfrac{N}{2})\Gamma(1-d+\tfrac{N}{2})}
{\Gamma(1-c+\tfrac{N}{2})\Gamma(b+\tfrac{N}{2})}\frac{\sin[\pi(d-\tfrac{N}{2})]}{\sin[\pi(c-\tfrac{N}{2})]}
= \Big(\frac{N}{2}\Big)^{\textup{i}(x+y)-2}(-1)^{N-m}(1+O(\tfrac{1}{N})),
\]
and, evidently, for $\arg(z)=\pi$ (i.e. $z=-N/2, \, N>0$) we have
$$
f\big(-\tfrac{N}{2}\big)=\Big(\frac{N}{2}\Big)^{\textup{i}(x+y)-2}(-1)^{N-n}(1+O(\tfrac{1}{N})).
$$

Using reflection properties \eqref{reflCgamma}, one can rewrite identity \eqref{CB} as
\beq
\int_{\mathbb{C}}[w-z_1]^{\alpha-1} [z_2-w]^{\beta-1} \frac{d^2w}{\pi}
= \frac{{\bf\tilde\Gamma}(\alpha,\beta,\gamma)}{[z_1-z_2]^{\gamma}},
\quad {\bf\tilde\Gamma}(\alpha_1,\ldots,\alpha_k):=\prod_{j=1}^k{\bf\Gamma}(\alpha_j|\alpha_j'),
\label{cbeta}\ee
where $\alpha+\beta+\gamma= \alpha'+\beta'+\gamma'=1$.
The integral converges for Re$(\alpha+\alpha')$, Re$(\beta+\beta')$, Re$(\gamma+\gamma')>0$.
Linear fractional transformations of the variables $w, z_1, z_2$ bring formula \eqref{cbeta}
to the following star-triangle relation form
\begin{eqnarray}
\int_{\mathbb{C}}[z_1-w]^{\alpha-1} [z_2-w]^{\beta-1} [z_3-w]^{\gamma-1}\frac{d^2w}{\pi}
=\frac{{\bf\tilde\Gamma}(\alpha,\beta,\gamma)} {[z_3-z_2]^{\alpha}[z_1-z_3]^{\beta}[z_2-z_1]^{\gamma}}.
 \label{STR}\end{eqnarray}
A multidimensional analogue of this formula --- a complex
 Selberg integral, was considered by Dotsenko and Fateev
in the context of two-dimensional conformal field theory \cite{DF} and by Aomoto  \cite{Aomoto}.

In \cite{DMV2017}, it was shown that  relation \eqref{STR} can be mapped to the following
infinite bilateral sum of  contour integrals (see below the Appendix for details)
\begin{eqnarray} && \makebox[-2em]{}
 \frac{1}{4\pi} \sum_{n\in \Z}
\int_{-\infty}^{\infty}\prod_{j=1}^3{\bf \Gamma}(b_j+y,n_j +n)
{\bf \Gamma}(a_j-y,n-m_j)dy  =\prod_{j,k=1}^3{\bf \Gamma}(b_j+a_k,n_j+m_k),
\label{MB}\end{eqnarray}
where Im$(b_j)$, Im$(a_j) < 0$ and the following balancing condition holds true
\begin{equation}
\sum_{j=1}^3(b_j+a_j)=-2\textup{i}, \qquad \sum_{j=1}^3(n_j+m_j)=0.
\label{balancing1}\end{equation}
A special case of this relation was formally obtained in \cite{BMS}.
The integrand in \eqref{MB} has poles at the following points of the integration variable
$$
y_p^{up}=-b_j + \textup{i}(|n+n_j|+2\Z_{\geq0}),\qquad y_p^{down}= a_j -\textup{i}(|n-n_j|+2\Z_{\geq0}),
$$
and the contour of integration should separate these sequences of poles
going upwards and downwards from the real line.
For completeness we mention that these are true poles remaining after partial cancellation
of zeros and poles of gamma functions forming the integrand,
so that there remains only the following set of integrand zeros
$$
y_z^{down}=-b_j + \textup{i}(n+n_j-2-2\Z_{\geq0}),\qquad y_z^{up}= a_j -\textup{i}(n-n_j-2-2\Z_{\geq0}).
$$
Asymptotic formulas given above for the products of complex gamma functions show
that the sum of integrals \eqref{MB} converges, since for $z=(n+\textup{i} y)/2\to\infty$
asymptotics of the integrand modulus is proportional to
$$
|z|^{\sum_{j=1}^3 [\textup{i}(b_j+a_j)-2]}= |z|^{-4}, \quad |z|\to\infty.
$$

Identity \eqref{MB} can be considered as a Mellin--Barnes form of the star-triangle relation \eqref{STR}, though
the character of counting independent free variables is different.
In \eqref{STR} one has two free complex variables (say, $z_1$ and $z_2$, since $z_3$ can be removed
by shifting the integration variable, $z_1$  and $z_2$), and two pairs of continuous and discrete
variables entering $\alpha$ and $\beta$. In \eqref{MB} one can shift the integration and summation
variables by arbitrary constants with an appropriate change of the integration contour, which
shows that one has actually only four independent pairs of the discrete and continuous variables.

In the following we consider multidimensional analogues of the exact formula \eqref{MB}
(as well as of a more general univariate relation derived in \cite{DM2019,GSVS1})
which can be considered as Mellin--Barnes forms of the multiple complex  hypergeometric integrals.
Also, we investigate more complicated higher order complex hypergeometric functions obeying
nontrivial symmetry transformations.

\section{Type I $C_n$ complex beta integral}

The univariate elliptic beta integral has been discovered in \cite{spi:umn}.
Its type I extension to root system $C_n$
was introduced by van Diejen and Spiridonov \cite{vDS2}. Different complete proofs
of its evaluation formula were given by Rains in \cite{rai:trans}
and by the second author in \cite{spi:short}. Therefore it can be called
as the van Diejen-Spiridonov-Rains $C_n$-integral of type I. It has the following structure.
Take the integration variables $z_j\in\mathbb{C}^\times,\, j=1,\ldots,n,$  and parameters
$t_1,\ldots,t_{2n+4}, p, q\in\mathbb{C}^\times$, satisfying restrictions $|p|, |q|, |t_j|<1$
and $\prod_{j=1}^{2n+4}t_j=pq$. Then one has the following exact integration formula
\begin{eqnarray}\nonumber &&
\kappa_n^C\int_{\T^n}\prod_{1\leq j<k\leq n}\frac{1}{\Gamma(z_j^{\pm 1} z_k^{\pm 1};p,q)}
\prod_{j=1}^n\frac{\prod_{\ell=1}^{2n+4}\Gamma(t_\ell z_j^{\pm 1};p,q)}
{\Gamma(z_j^{\pm2};p,q)}\, \frac{dz_j}{z_j}
\\ && \makebox[4em]{}
=\prod_{1\leq \ell<s\leq 2n+4}\Gamma(t_\ell t_s;p,q),
\label{C-typeI}\end{eqnarray}
where $\T$ is the unit circle of positive orientation,
$$
\kappa_n^C=\frac{(p;p)_\infty^n(q;q)_\infty^n}{(4\pi \textup{i})^n n!},
\qquad (t;q)_\infty=\prod_{k=0}^\infty (1-tq^k),
$$
and
$$
\Gamma(z;p,q):=
\prod_{j,k=0}^\infty\frac{1-z^{-1}p^{j+1}q^{k+1}}{1-zp^{j}q^{k}}
$$
is the elliptic gamma function with the standard convention
\begin{eqnarray*} &&
\Gamma(tz^{\pm1};p,q):=\Gamma(tz;p,q)\Gamma(tz^{-1};p,q),
\\ &&
\Gamma(tz_j^{\pm1}z_k^{\pm1};p,q):=\Gamma(tz_jz_k;p,q)\Gamma(tz_jz_k^{-1};p,q)
\Gamma(tz_j^{-1}z_k;p,q)\Gamma(tz_j^{-1}z_k^{-1};p,q).
\end{eqnarray*}

The following asymptotic relation for the elliptic gamma function holds true \cite{ruij}
(see also \cite{vDS3})
\beq\label{parlim2}
\Gamma(e^{-2\pi vu};e^{-2\pi v\omega_1},e^{-2\pi v\omega_2})
\stackreb{\sim}{ v\to 0^+}
e^{-\pi\frac{2u-\omega_1-\omega_2}{12v\omega_1\omega_2}}\gamma^{(2)}(u;\omega_1,\omega_2),
\ee
where $ \gamma^{(2)}(u;\omega_1,\omega_2)$ is the Faddeev modular quantum dilogarithm \cite{Fad94,Fad95}
called also the hyperbolic gamma function \cite{ruij}
\beq
\gamma^{(2)}(u;\mathbf{\omega})= \gamma^{(2)}(u;\omega_1,\omega_2):=e^{-\frac{\pi\textup{i}}{2}
B_{2,2}(u;\mathbf{\omega}) } \gamma(u;\mathbf{\omega}),
\label{HGF}\ee
with the multiple Bernoulli polynomial of the  second order
$$
 B_{2,2}(u;\mathbf{\omega})=\frac{1}{\omega_1\omega_2}
\left((u-\frac{\omega_1+\omega_2}{2})^2-\frac{\omega_1^2+\omega_2^2}{12}\right)
$$
and
\beq
\gamma(u;\mathbf{\omega}):= \frac{(\tilde {\bf q} e^{2\pi \textup{i} \frac{u}{\omega_1}};\tilde {\bf q})_\infty}
{(e^{2\pi \textup{i} \frac{u}{\omega_2}};{\bf q})_\infty}
=\exp\left(-\int_{\R+\textup{i}0}\frac{e^{ux}}
{(1-e^{\omega_1 x})(1-e^{\omega_2 x})}\frac{dx}{x}\right).
\label{int_rep}\ee
For {Im}$(\omega_1/\omega_2)>0$, one has the parametrization
\beq
{\bf q}= e^{2\pi \textup{i}\frac{\omega_1}{\omega_2}},\qquad \tilde {\bf q}= e^{-2\pi \textup{i}\frac{\omega_2}{\omega_1}}.
\label{bfq}\ee
This $\bf q$ is different from $q$ figuring in the definition of elliptic hypergeometric integral \eqref{C-typeI}.
For $|{\bf q}|<1$ the infinite product representation defines $\gamma^{(2)}(u;\mathbf{\omega})$
as a meromorphic function of $u\in\CC$.
For $|{\bf q}|=1$ one should use the integral representation converging in a restricted domain of $u$
(say, for $\omega_1, \omega_2>0$ one should have $0<\text{Re}(u)< \omega_1+\omega_2$).
Note also that though the parametrization of $\bf q$ \eqref{bfq} is not symmetric in
$\omega_1$ and $\omega_2$, the resulting function \eqref{HGF} obeys this symmetry.

As shown by Rains in \cite{rai:limits}, the limit \eqref{parlim2} is
uniform on compacta with exponentially small corrections. Therefore applying
the following parametrization of variables
$$
t_\ell=e^{-2\pi v g_\ell}, \qquad z_j=e^{-2\pi v u_j},\qquad p=e^{-2\pi v\omega_1}, \quad q=e^{-2\pi v\omega_2},
\; v>0
$$
(which assumes that Re$(\omega_1)$, Re$(\omega_2) > 0$), one can establish rigorously the asymptotic
limit $v\to 0^+$ for all known elliptic hypergeometric integrals.
In particular, from \eqref{C-typeI} one obtains in this way the following identity
\begin{eqnarray} \nonumber &&
\frac{1}{n!}\int_{u_j\in\textup{i}\R}\prod_{1\leq j<k\leq n}\frac{1}{\gamma^{(2)}(\pm u_j\pm u_k;\mathbf{\omega})}
\prod_{j=1}^n\frac{\prod_{\ell=1}^{2n+4}\gamma^{(2)}(g_\ell\pm u_j;\mathbf{\omega}) }
{\gamma^{(2)}(\pm 2 u_j;\mathbf{\omega}) } \frac{du_j}{2\textup{i}\sqrt{\omega_1\omega_2}}
\\  && \makebox[4em]{}
=\prod_{1\leq \ell<s\leq 2n+4}\gamma^{(2)}(g_\ell+g_s;\mathbf{\omega}),
\label{hyperCI}\end{eqnarray}
where  the following balancing condition holds true
\beq\label{balcon}
\sum_{\ell=1}^{2n+4} g_\ell=Q, \qquad Q:=\omega_1+\omega_2.
\ee
Here we use the compact notation $\gamma^{(2)}(g\pm u;\mathbf{\omega})
:=\gamma^{(2)}(g+u;\mathbf{\omega})\gamma^{(2)}(g-u;\mathbf{\omega})$.
The contour of integration in \eqref{hyperCI} can be chosen as the imaginary axis
under the conditions Re$(g_\ell)>0$ that stem from the requirement
$|t_\ell|<1$ at the elliptic level.
For further considerations it is convenient to assume that $\omega_1\omega_2\in\R_{> 0}$,
which can always be realized using the homogeneity $\gamma^{(2)}(\lambda u;\lambda\omega_1,\lambda\omega_2)
=\gamma^{(2)}(u;\omega_1,\omega_2)$ for $\lambda\neq 0$.
In the conformal field theory the canonical normalization is $\omega_1\omega_2=1$.

Since the poles of $\gamma^{(2)}(u;\mathbf{\omega})$ are located at
$u_p\in \{ -n\omega_1 -m\omega_2\},\, n,m \in \Z_{\geq0},$
the poles of the integrand in \eqref{hyperCI}  are
$$
u_j^{poles}\in \{ g_\ell+\omega_1\Z_{\geq0}+ \omega_2\Z_{\geq0}\}\cup
\{ -g_\ell-\omega_1\Z_{\geq0}- \omega_2\Z_{\geq0}\},
\quad \ell=1,\ldots, 2n+4.
$$
The contours of integration in \eqref{hyperCI} should separate these two sets of points
and, therefore, they differ from the imaginary axis, if the conditions Re$(g_\ell)>0$
are violated.

The function \eqref{int_rep} plays an important role in the analysis
of $2d$  quantum Liouville theory. In the corresponding literature one denotes
$b:=\sqrt{\frac{\omega_1}{\omega_2}}$, i.e. ${\bf q}=e^{2\pi \textup{i} b^2}$
and $\tilde {\bf q}=e^{-2\pi \textup{i} b^{-2}}$, so that the central charge $c$ of the
theory is parametrized as $c=1+6(b+b^{-1})^2$ \cite{DF}.

In \cite{GSVS1}, the simultaneous limit ${\bf q}\to 1$ and $\tilde {\bf q}\to 1$
(when $b\to \textup{i}$, or $c\to 1$) was rigorously investigated and corresponding degeneration of
some particular hyperbolic integrals was considered. Take small $\delta>0$  and set  $b=\textup{i}+\delta$.
Then  for $\delta\to 0^+$ the following asymptotic limit holds true
\beq\label{gam2lim2}
\gamma^{(2)}(\textup{i}\sqrt{\omega_1\omega_2}(n+x\delta);\omega_1,\omega_2)\stackreb{\sim}{\delta\to 0^+} e^{\frac{\pi \textup{i}}{2}n^2} (4\pi\delta)^{\textup{i}x-1}{\bf \Gamma}(x,n),
\quad \sqrt{\omega_1\over \omega_2}=\textup{i}+\delta,
\ee
where $n\in \Z, \, x\in\CC$, and ${\bf \Gamma}(x,n)$ is the complex gamma function \eqref{Cgamma}.
This limit is uniform on compacta with exponentially small corrections. Therefore all
degenerations of the hyperbolic integrals considered below rigorously prove the
resulting limiting relations.

Now we apply this limit to the multiple beta integral \eqref{hyperCI},
which was done earlier in the univariate setting in \cite{GSVS2}.
The structure of poles $u_j^{poles}$ shows that if parameters $g_\ell$
approach the values $\textup{i}\sqrt{\omega_1\omega_2}N_\ell$ with $N_\ell$ integers
or half-integers, then for small $\delta$ infinitely many poles start to pinch the integration
contours at the points $u_j=\textup{i}\sqrt{\omega_1\omega_2}\, N$, $N\in \Z+\nu,\, \nu =0 $
or $\tfrac12$, respectively.
As demonstrated in \cite{GSVS2}, for a function $\Delta(u)$ formed from the $\gamma^{(2)}(g+u;\mathbf{\omega})$-functions with different parameters $g$,
the asymptotics of its integral over the imaginary axis in such a regime has the form
\beq\label{intdel}
\int_{-\textup{i}\infty}^{\textup{i}\infty}\Delta(u){du\over \textup{i}\sqrt{\omega_1\omega_2}}
\stackreb{\sim}{\delta\to 0^+}
\sum_{m\in \Z+\nu}
\int_{-\infty}^{\infty} \left[
\stackreb{\lim}{\delta\to 0}\delta\Delta(\textup{i}\sqrt{\omega_1\omega_2}(m+x\delta))\right] dx,
\ee
provided on the right-hand side one has a well defined asymptotic limit. The symbol $\sim$ in
\eqref{intdel} means that the ratio of the left- and right-hand side expressions goes to 1 for
$\delta\to 0^+$. For uniformness
of the limit there should emerge a power of $\delta$ multiplied by a converging infinite
bilateral sum of Mellin--Barnes type integrals. For multiple integrals, one has the
asymptotic relation \eqref{intdel} for each integration variable.

Now we set $\sqrt{\frac{\omega_1}{\omega_2}}=\textup{i}+\delta$ in  \eqref{hyperCI},
parametrize the integration variables as
\beq\label{zkom}
u_j= \textup{i}\sqrt{\omega_1\omega_2}(m_j+\delta x_j),
\quad x_j\in\CC,\quad m_j\in\Z+\nu,\; \nu=0, \frac{1}{2},
\ee
the constants $g_\ell$ as
\beq\label{gkom}
g_\ell=\textup{i}\sqrt{\omega_1\omega_2}(N_\ell+\delta a_\ell),\quad a_\ell\in\CC,
\quad N_\ell\in\Z+\nu,\; \nu=0, \frac{1}{2},
\ee
and take the limit $\delta\to 0^+$.
 Then the balancing condition $\sum_{k=1}^{2n+4} g_k=\omega_1+\omega_2$
splits into two separate constraints on $a_k$ and $N_k$:
\beq\label{xknk}
\sum_{k=1}^{2n+4}a_k=-2\textup{i},\qquad \sum_{k=1}^{2n+4} N_k=0.
\ee
Since we took $\omega_1\omega_2>0$, the conditions Re$(g_\ell)>0$ pass to the constraints Im$(a_\ell)<0$.
The discrete parameter $\nu=0, \frac{1}{2}$  emerged from the demand for the combinations
$m_j\pm N_k$ to take integer values, as required in  \eqref{gam2lim2}.

We recall that the function $\gamma^{(2)}(y;\omega_1,\omega_2)$ has the following asymptotics:
\begin{eqnarray} \nonumber &&
{\rm lim}_{y\to \infty}\gamma^{(2)}(y;\omega_1,\omega_2)=e^{-{\textup{i} \pi\over 2}B_{2,2}(y,\omega_1,\omega_2)},
\quad {\rm for}\; {\rm arg}\;\omega_1<{\rm arg}\; y<{\rm arg}\;\omega_2+\pi,
\\ &&
{\rm lim}_{y\to \infty}\gamma^{(2)}(y;\omega_1,\omega_2)=e^{{\textup{i} \pi\over 2}B_{2,2}(y,\omega_1,\omega_2)},
\quad {\rm for}\; {\rm arg}\;\omega_1-\pi<{\rm arg}\; y<{\rm arg}\;\omega_2.
\label{asy1}\end{eqnarray}
It follows that the integrand of the original hyperbolic integral for $x_j\to \pm \infty$
has the asymptotics $e^{- 12\pi \delta \sum_{j=1}^n |x_j|}$. For $\delta\to 0^+$ the integral diverges and
our goal is to estimate the rate of this divergence.

Inserting parametrizations (\ref{zkom}) and (\ref{gkom}) in (\ref{hyperCI}), and using
formula (\ref{gam2lim2}), we find the limiting relations
$$
\prod_{\ell=1}^{2n+4}\gamma^{(2)}(g_\ell\pm u_j;\mathbf{\omega})\to \frac{(-1)^{2\nu n}}{(4\pi\delta)^{4(n+1)}}
\prod_{\ell=1}^{2n+4}{\bf \Gamma}(a_\ell\pm x_j,N_\ell\pm m_j),
$$
$$
\prod_{1\leq \ell < s\leq 2n+4}
\gamma^{(2)}(g_\ell+g_s;\mathbf{\omega})\to
\frac{(-1)^{(n+1)(n+2)\nu}}{(4\pi\delta)^{n(2n+3)}}\prod_{1\leq \ell < s \leq 2n+4}
{\bf \Gamma}( a_\ell+a_s,N_\ell+N_s),
$$
$$
\prod_{1\leq j<k\leq n}\gamma^{(2)}(\pm u_j\pm u_k;\mathbf{\omega}) \to
\frac{2^{2n(n-1)}}{(4\pi\delta)^{2n(n-1)}}
\frac{1}{\prod_{1\leq j<k\leq n} [(x_j\pm x_k)^2+(m_j\pm m_k)^2]},
$$
\beq
\gamma^{(2)}(\pm 2u_j;\mathbf{\omega}) \to \frac{(-1)^{2\nu}}{(4\pi\delta)^{2}}
{\Gamma(m_j+ \textup{i}x_j)\over \Gamma(1+m_j- \textup{i}x_j)}
{\Gamma(-m_j- \textup{i}x_j)\over \Gamma(1-m_j+ \textup{i}x_j)}={(4\pi\delta)^{-2}\over x_j^2+m_j^2}.
\label{extrasign}\ee
Collecting all the multipliers and cancelling
$(-1)^{2n\nu}(4\pi\delta)^{-n(2n+3)}$ --- the common diverging factor
on both sides of the equality  \eqref{hyperCI}, we obtain our key complex beta integral:
\begin{eqnarray} \label{complexCI} &&
\frac{1}{(2^{2n+1}\pi)^n n!}\sum_{m_j\in \Z+\nu}\int_{x_j\in\R}
\prod_{1\leq j<k\leq n}\left[(x_j\pm x_k)^2+(m_j\pm m_k)^2\right]
\prod_{j=1}^n (x_j^2+m_j^2)
\\  && \makebox[-2em]{} \times
\prod_{j=1}^n\prod_{\ell=1}^{2n+4}{\bf \Gamma}(a_\ell\pm x_j,N_\ell\pm m_j)dx_j
=(-1)^{(n-1)(n-2)\nu}\prod_{1\leq \ell<s\leq 2n+4}{\bf \Gamma}(a_\ell+a_s,N_\ell+N_s),\quad
\nonumber \end{eqnarray}
where $\sum_{\ell=1}^{2n+4} a_\ell=-2\textup{i}$, $\sum_{\ell=1}^{2n+4} N_\ell=0$, and
$$
{\bf \Gamma}(x_1\pm x_2 ,n_1\pm n_2):={\bf \Gamma}(x_1+ x_2 ,n_1+ n_2){\bf \Gamma}(x_1- x_2, n_1- n_2).
$$
Here we have the variables $N_\ell, m_j \in\Z+\nu,\, \nu=0,\,\frac{1}{2}$, so that
$N_\ell\pm m_j$ and $N_\ell+N_s$  take integer values. This formula was first obtained
by a different method in \cite{DM2019} (see formula (3.7) there).

Let us check  that the infinite bilateral sums of integrals \eqref{complexCI} converge.
According to \eqref{asympgam}, when $z_j=\tfrac12(m_j+\textup{i}x_j )\to\infty$, $j$ fixed,
the modulus of the integrand behaves as
$$
|z_j|^{\sum_{\ell=1}^{2n+4} (2\textup{i} a_\ell-2)+4(n-1)+2}=|z_j|^{-6}
$$
due to the balancing condition. If simultaneously several $z_j$-variables go to infinity,
then the integrand vanishes even faster, since moduli of the cross-terms $(x_j\pm x_k)^2+(m_j\pm m_k)^2$
grow slower than $|z_jz_k|^{2}$. This guarantees convergence of the sum of integrals.

For any $m_j$ true poles of the integrand in \eqref{complexCI} are located at the points
\beq\label{poles}
x_j^{(p)} \in \{ \textup{i}|m_j+N_k|-a_k+2\textup{i}\Z_{\geq0},\} \cup
\{-\textup{i}|m_j-N_k|+a_k-2\textup{i}\Z_{\geq0}\}.
\ee
The contours of integration over $x_j$ in \eqref{complexCI} should separate these sets
of poles. The real axes are valid for Im$(a_k)<0$ (the conditions following from the
demand $\text{Re}(g_k/\sqrt{\omega_1\omega_2})>0$ and $\delta\to 0^+$).
As to the cases $a_j\in\R$, $a_j\neq - a_k$, we can perform analytical continuation by deforming
the contour of integration slightly below and above the real axis at the appropriate points
in such a way that for Im$(a_k)=0$ no poles emerge on the integration contour, and the formula
remains true.

Formula \eqref{complexCI} represents a complex
analogue of the Gustafson $C_n$-integral evaluation formula described in Theorem 5.3
of \cite{gus:some1}. However, it has a substantially more symmetric  form due to the
balancing condition symmetric in all parameters, which was not possible to realize
in Gustafson's case. Note that
the original elliptic hypergeometric integral \eqref{C-typeI} can be reduced to the
Gustafson's one by taking two different limits. In the first case, one should express
$t_{2n+4}$ in terms of other parameters with the help of the balancing condition, apply
the elliptic gamma function reflection formula to remove $p$-dependence in the
arguments, and take the $p\to 0$ limit with all remaining free parameters $t_j$
being fixed. After that, one should take the limit $q\to 1$ in the arising
trigonometric $q$-hypergeometric integral with $t_j=q^{a_j}$ for fixed $a_j$, which ends
up in the relation of interest. In the second case, one goes first from the elliptic
to hyperbolic level as described above. Then it is necessary to take the limit $b\to0$
(when ${\bf q}\to 1$ and  $\tilde {\bf q}\to 0$) after
expressing one of the parameters in terms of others with the help of the balancing condition
and applying the reflection formula for the
$\gamma^{(2)}$-function. The latter possibility was described in \cite{rai:limits}
or in our recent paper \cite{GSVS3}.

\section{Type I $A_n$ complex beta integral}

A type I multiple elliptic beta integral on the root system $A_n$ was suggested by the second author
in \cite{spi:theta2}. Different complete proofs of its evaluation formula were given by Rains \cite{rai:trans}
and Spiridonov in \cite{spi:short}. Therefore it can be called the Rains-Spiridonov $A_n$-integral of type I.
To describe it, we take parameters $p, q, t_\ell, s_\ell\in\mathbb{C}^\times,\, \ell=1,\ldots, n+2,$
satisfying the constraints $|p|, |q|, |t_m|, |s_m|<1$ and $ST=pq$, where
$S=\prod_{\ell=1}^{n+2}s_\ell$, $T=\prod_{\ell=1}^{n+2}t_\ell$. Then,
\begin{eqnarray}\nonumber
&& \kappa_{n}^A\int_{{\mathbb T}^n}
\prod_{1\le j<k\le n+1}\frac{1}{\Gamma(z_jz_k^{-1},z_j^{-1}z_k;p,q)}
\,\prod_{j=1}^{n+1}\prod_{\ell=1}^{n+2}\Gamma(s_\ell z_j,t_\ell z_j^{-1};p,q)
\,\prod_{i=1}^n\frac{dz_i}{z_i}
\\ && \makebox[4em]{}
=\prod_{\ell=1}^{n+2} \Gamma(Ss_\ell^{-1},Tt_\ell^{-1};p,q)
\prod_{\ell,s=1}^{n+2} \Gamma(s_\ell t_s;p,q),
\label{AI}\end{eqnarray}
where $z_1z_2\cdots z_{n+1}=1$ and
$$
\kappa_{n}^A=\frac{(p;p)_\infty^n(q;q)_\infty^n}{(2\pi \textup{i})^n(n+1)!}.
$$
Under the taken restrictions on the parameters the integral kernel in \eqref{AI} has
double sequences of poles in each integration variable $z_j=t_\ell p^aq^b,\, j=1,\ldots, n,$
and $z_1\cdots z_n= s_\ell p^aq^b,\, a,b\in\Z_{\geq 0}$, which accumulate near zero.
Since for a fixed $z_j$ other integration variables $z_k,\, k \neq j,$ lie on the unit circle,
these poles lie inside $\T$. In a similar way, other sequences of poles
$z_j=t_\ell^{-1} p^{-a}q^{-b},\, j=1,\ldots, n,$ and
$z_1\cdots z_n= s_\ell^{-1} p^{-a}q^{-b}$, lie outside $\T$ and go to infinity as $a, b \to\infty$.

Similarly to the previous case, reduction to the hyperbolic level yields the identity
\begin{eqnarray} \nonumber &&
\frac{1}{(n+1)!}\int_{u_j\in\textup{i}\R}
\frac{\prod_{j=1}^{n+1}\prod_{\ell=1}^{n+2}\gamma^{(2)}(g_\ell+ u_j,f_\ell-u_j;\mathbf{\omega}) }
{\prod_{1\leq j<k\leq n+1}\gamma^{(2)}(\pm(u_j- u_k);\mathbf{\omega})
 } \prod_{j=1}^n \frac{du_j}{\textup{i}\sqrt{\omega_1\omega_2}}
\\  && \makebox[4em]{}
=\prod_{\ell=1}^{n+2}\gamma^{(2)}(G-g_\ell, F-f_\ell;\mathbf{\omega})
\prod_{\ell,s=1}^{n+2}\gamma^{(2)}(g_\ell+f_s;\mathbf{\omega}),
\label{hyperAI}\end{eqnarray}
where Re$(\omega_{1,2})$, Re$(g_\ell)$, Re$(f_\ell)>0$, $\sum_{j=1}^{n+1}u_j=0$, $G=\sum_{\ell=1}^{n+2}g_\ell$ and $F=\sum_{\ell=1}^{n+2}f_\ell$ with $G+F=\omega_1+\omega_2$. Integrations over $u_j$ go along
the imaginary axes which separate poles going to infinity to the right,
$u_j= f_\ell +\omega_1\Z_{\geq0} +\omega_2\Z_{\geq0}$
and $\sum_{k=1}^n u_k= g_\ell +\omega_1\Z_{\geq0} +\omega_2\Z_{\geq0},$ from the poles
$u_j= -g_\ell -\omega_1\Z_{\geq0} -\omega_2\Z_{\geq0},\, \sum_{k=1}^nu_k
= -f_\ell -\omega_1\Z_{\geq0} -\omega_2\Z_{\geq0},$ going to infinity to the left.

For the transition to the level of complex hypergeometric functions along the lines
indicated in the previous section, we parametrize the integration variables as
\beq\label{zkomA}
u_j= \textup{i}\sqrt{\omega_1\omega_2}(m_j+\delta x_j),
\quad x_j\in\CC,\quad m_j\in\Z+\nu,\; \nu=0, \frac{1}{2},
\ee
and the constants $g_\ell, f_\ell$ as
\beq\label{gkomA}
g_\ell=\textup{i}\sqrt{\omega_1\omega_2}(N_\ell+\delta a_\ell),\qquad
f_\ell=\textup{i}\sqrt{\omega_1\omega_2}(M_\ell+\delta b_\ell),
\ee
where $a_\ell, b_\ell \in\CC$ and $N_\ell, M_\ell\in\Z+\nu,\; \nu=0, \frac{1}{2},$
and take the limit $\delta\to 0^+$. Then the parameters $a_\ell, b_\ell, N_\ell, M_\ell$
satisfy the constraints Im$(a_\ell)$, Im$(b_\ell)<0$ (stemming from the original restrictions
on $g_\ell$ and $f_\ell$) as well as the balancing conditions
\beq
A + B=-2\textup{i},\quad A=\sum_{\ell=1}^{n+2}a_\ell,\quad B=\sum_{\ell=1}^{n+2}b_\ell,
\label{xknkAA}\ee
and
\beq
N+M=0,\quad N=\sum_{\ell=1}^{n+2} N_\ell,\quad M=\sum_{\ell=1}^{n+2} M_\ell,
\label{xknkA}\ee
following from the original constraint $G+F=\omega_1+\omega_2$.

Inserting expressions (\ref{zkomA}) and (\ref{gkomA}) in (\ref{hyperAI}), and recalling
 the asymptotics (\ref{gam2lim2}), we find the limiting relations
\begin{eqnarray*} &&
\prod_{j=1}^{n+1}\prod_{\ell=1}^{n+2}\gamma^{(2)}(g_\ell+ u_j,f_\ell-u_j;\mathbf{\omega})\to
\frac{e^{\pi\textup{i}(n+1)(\frac{1}{2}\sum_{\ell=1}^{n+2}(N_\ell^2+M_\ell^2)+(n+2)\nu^2)}}
{(4\pi\delta)^{2(n+1)^2}}
\\ && \makebox[2em]{} \times
\prod_{j=1}^{n+1}\prod_{\ell=1}^{n+2}{\bf \Gamma}(a_\ell+ x_j,N_\ell+m_j){\bf \Gamma}(b_\ell- x_j,M_\ell-m_j),
\\ && \makebox[2em]{}
\prod_{1\leq j<k\leq n+1}\gamma^{(2)}(\pm (u_j- u_k);\mathbf{\omega}) \to
\frac{2^{n(n+1)}(4\pi\delta)^{-n(n+1)}}
{\prod_{1\leq j<k\leq n+1} [(x_j-x_k)^2+(m_j-m_k)^2]},
\\ &&
\prod_{\ell=1}^{n+2}\gamma^{(2)}(G-g_\ell, F-f_\ell;\mathbf{\omega})
\prod_{\ell,s=1}^{n+2}\gamma^{(2)}(g_\ell+f_s;\mathbf{\omega})\to
{e^{\textup{i}{\pi\over 2}( 2N^2(n-1)+(n+3)\sum_{\ell=1}^{n+2}(N_\ell^2+M_\ell^2))}\over (4\pi\delta)^{n^2+2n+2}}
\\ \nonumber  && \makebox[2em]{} \times
\prod_{\ell=1}^{n+2}{\bf \Gamma}(A-a_\ell,N-n_\ell){\bf \Gamma}( B-b_\ell,M-M_\ell)
\prod_{\ell,s=1}^{n+2}{\bf \Gamma}(a_\ell+b_s,N_\ell+M_s).
\end{eqnarray*}
Collecting all the multipliers and cancelling
 the common diverging factor
on both sides of the equality  \eqref{hyperAI}, we obtain
in the limit $\delta\to 0^+$   the type I complex beta integral for the root system $A_n$:
\begin{eqnarray} \nonumber &&
\frac{1}{(2^{n+3}\pi)^n (n+1)!}\sum_{m_j\in \Z+\nu}\int_{x_j\in\R}
\prod_{j=1}^{n+1} \prod_{\ell=1}^{n+2}{\bf \Gamma}(a_\ell+ x_j,N_\ell+ m_j){\bf \Gamma}(b_\ell- x_j,M_\ell- m_j)
\\ \nonumber && \makebox[1em]{}  \times
\prod_{1\leq j<k\leq n+1}\left[(x_j-x_k)^2+(m_j-m_k)^2\right] \prod_{j=1}^n dx_j
=e^{\textup{i} \pi(n-1)(N^2-(n+2)\nu^2)}
\\  && \makebox[2em]{} \times
\prod_{\ell=1}^{n+2}{\bf \Gamma}(A-a_\ell,N-N_\ell){\bf \Gamma}( B-b_\ell,M-M_\ell)
\prod_{\ell,s=1}^{n+2}{\bf \Gamma}(a_\ell+b_s,N_\ell+M_s),
\label{complexAI}\end{eqnarray}
where $\sum_{j=1}^{n+1}x_j=\sum_{j=1}^{n+1}m_j=0$.
For even dimensional integrals this relation holds true only for $\nu=0$.
However, for odd $n$ both values $\nu=0 $ and $\nu=1/2$ are admissible.
This new complex hypergeometric beta integral evaluation formula was announced in \cite{GSVS4}.

The integrations in \eqref{complexAI} go along the real axes which separate
the poles $x_j= b_\ell -\textup{i}(|m_j-M_\ell|+2\Z_{\geq0}), \, j=1,\ldots,n,$
and $\sum_{k=1}^n x_k=a_\ell -\textup{i}(|\sum_{k=1}^n m_k-N_\ell|+2\Z_{\geq0})$ going to infinity
downwards the real line from the poles $x_j= -a_\ell +\textup{i}(|m_j+N_\ell|+2\Z_{\geq0}), \, j=1,\ldots,n,$
and $\sum_{k=1}^n x_k=-b_\ell  +\textup{i}(|\sum_{k=1}^n m_k+M_\ell|+2\Z_{\geq0})$ going to infinity upwards.

The asymptotic behavior of the integrand modulus in \eqref{complexAI} for
$z_j=\tfrac12(m_j+\textup{i}x_j )\to\infty$, $j\leq n$ fixed, can be found using \eqref{asympgam}.
Since $z_{n+1}=-\sum_{j=1}^{n}z_j$ we have contributions of two kinds, from the terms containing only $z_j$ and terms containing $z_{n+1}$:
$$
\propto|z_j|^{\sum_{\ell=1}^{n+2} (\textup{i}( a_\ell+b_\ell)-2)+2(n-1)}=|z_j|^{-4}
\quad \text{and} \quad
\propto|z_{n+1}|^{\sum_{\ell=1}^{n+2} (\textup{i}( a_\ell+b_\ell)-2)+2n}=|z_{n+1}|^{-2},
$$
where we used the balancing condition \eqref{xknkAA}. So, for large values of all $|z_k|$ the integrand
behaves as least as $\prod_{k=1}^n |z_k|^{-6}$ in integration variables and the infinite
sums of integrals of interest converge.

Formula \eqref{complexAI}  represents a complex analogue of the Gustafson $A_n$-integral evaluation
described in Theorem 5.1 of \cite{gus:some1}, although it has substantially more symmetric
form due to the balancing condition symmetric in all parameters, which was not possible
to have in Gustafson's case.
In the same way as in the $C_n$-integrals case, Gustafson's $A_n$-integral evaluation formula
can be obtained in two different limits from the original elliptic hypergeometric integral
which we shall not describe here, since the difference carries only a technical character.

\section{Type II integrals on the root system $A_n$}

Recently, in \cite{GSVS3} we derived a generalized complex Selberg integral in the
Mellin--Barnes form. It was obtained in a chain of reductions of the type II elliptic
beta integral on the root system $C_n$ evaluated by van Diejen and Spiridonov \cite{vDS2}
as a consequence of relation \eqref{C-typeI}. Explicitly, it has the following form
\begin{eqnarray} \nonumber &&
\frac{1}{(8\pi)^n n!} \sum_{m_j\in \Z+\mu}\int_{u_j\in \R}
\prod_{1\le j<k\le n}\! \frac{{\bf \Gamma}(\gamma \pm u_j\pm u_k, r\pm m_j\pm m_k)}
{{\bf \Gamma}(\pm u_j\pm u_k, \pm m_j\pm m_k)}
\\ && \makebox[6em]{} \times
\ \prod_{j=1}^n  \left[ \prod_{\ell=1}^6{\bf \Gamma}(\gamma_\ell\pm u_j, r_\ell\pm m_j)\right]
(u_j^2+m_j^2)\, du_j
\label{complexS} \\ && \makebox[0em]{}
= (-1)^{r\frac{n(n-1)}{2})} \prod_{j=1}^n\frac{{\bf \Gamma}(j\gamma,jr)}{{\bf \Gamma}(\gamma,r)}
\prod_{1\leq \ell<s\leq 6}{\bf \Gamma}((j-1)\gamma+\gamma_\ell+\gamma_s,(j-1)r+r_\ell+r_s),
\nonumber \end{eqnarray}
where $r\in\Z,\,  m_j, r_\ell\in \Z+\mu,\, \mu=0, \frac{1}{2},$ and the continuous
variables $\gamma, \gamma_\ell\in\CC$ satisfy the balancing condition
\begin{equation}
(2n-2)\gamma+\sum_{\ell=1}^6\gamma_\ell=-2\textup{i}, \quad (2n-2)r+\sum_{\ell=1}^6r_\ell=0.
\label{bal2}\end{equation}
Here we use the compact notation
\begin{eqnarray*} &&
{\bf \Gamma}(\gamma \pm a\pm b, r\pm m\pm n):={\bf \Gamma}(\gamma + a+ b, r+ m+ n)
{\bf \Gamma}(\gamma + a- b, r+ m- n)
\\ && \makebox[3em]{} \times
{\bf \Gamma}(\gamma - a+ b, r- m+ n)
{\bf \Gamma}(\gamma - a-b, r- m- n).
\end{eqnarray*}
Similar to the previous cases, large $|z_j|=\tfrac12|m_j+\textup{i} u_j|$ behavior of the integrand modulus
in \eqref{complexS} is $\propto |z_j|^{4i\gamma(n-1)+2i\sum_1^6\gamma_l-12+2}=|z_j|^{-6}$
 and the integrals, as well as their infinite sums converge.

In this section we discuss the type II elliptic beta integrals on the root system $A_n$
that were constructed in \cite{spi:theta2}. As shown in \cite{SV} these integrals represent a
generalization of the elliptic Selberg integral introduced in \cite{vDS1,vDS2}. Since these integrals do
not appear anywhere else than in \cite{spi:theta2} and some surveys \cite{spi:essays,SV},
they can be called Spiridonov's type II integrals on the root system $A_n$.

Let us describe these integrals explicitly and perform their reduction to the
level of complex hypergeometric functions (such a reduction for the elliptic
Selberg integral was performed in \cite{GSVS3}).
We take 10 parameters $p,q,$ $t, s,$ $t_j,$ $s_j$, $j=1,2,3,$
such that $|p|, |q|, |t_j|, s_j|<1$ and $(ts)^{n-1}\prod_{k=1}^3t_ks_k=pq$ (the balancing condition).
Then the $A_n$-integrals of interest have the form
\begin{eqnarray} \nonumber && \makebox[-2em]{}
I_n(\underline{t},\underline{s};p,q)=
\kappa_{n}^A \int_{\T^n} \prod_{1\leq i<j\leq n+1}
\frac{\Gamma(tz_iz_j,sz_i^{-1}z_j^{-1};p,q)}
{\Gamma(z_iz_j^{-1},z_i^{-1}z_j;p,q)}
\\  && \makebox[4em]{} \times
\prod_{j=1}^{n+1}\prod_{k=1}^3\Gamma(t_kz_j,s_kz_j^{-1};p,q)
\prod_{j=1}^n \frac{dz_j}{z_j},
 \label{AIIa} \end{eqnarray}
where $\prod_{j=1}^{n+1}z_j=1$. These integrals can be computed explicitly.
For odd $n$ one has
\begin{eqnarray} \label{AIIa-odd} &&
I_n(\underline{t},\underline{s};p,q)=
\Gamma(t^{\frac{n+1}{2}},s^{\frac{n+1}{2}};p,q)
\prod_{1\leq i<k\leq 3} \Gamma(t^{\frac{n-1}{2}}t_it_k,
s^{\frac{n-1}{2}}s_is_k;p,q)
\\ && \makebox[-1em]{}  \times
\prod_{\ell=1}^{(n+1)/2}\prod_{i,k=1}^3\Gamma((ts)^{\ell-1}t_is_k;p,q)
\prod_{j=1}^{(n-1)/2}\Gamma((ts)^j;p,q)\prod_{1\leq i<k\leq 3}
\Gamma(t^{j-1}s^jt_it_k,t^js^{j-1}s_is_k;p,q),
 \nonumber\end{eqnarray}
and for even $n$ one has
\begin{eqnarray}\label{AIIa-even} &&
I_n(\underline{t},\underline{s};p,q)
= \prod_{i=1}^3\Gamma(t^{\frac{n}{2}}t_i,s^{\frac{n}{2}}s_i;p,q)
\Gamma(t^{\frac{n}{2}-1}t_1t_2t_3,s^{\frac{n}{2}-1}s_1s_2s_3;p,q)
\\ && \makebox[-1em]{} \times
\prod_{j=1}^{n/2}\Gamma((ts)^j;p,q)
\prod_{i,k=1}^3\Gamma((ts)^{j-1}t_is_k;p,q)
\prod_{1\leq i<k\leq 3}\Gamma(t^{j-1}s^jt_it_k,t^js^{j-1}s_is_k;p,q).
\nonumber\end{eqnarray}

The hyperbolic reduction of integrals \eqref{AIIa}, performed
completely similarly to the previous cases, has the following form
\begin{eqnarray} \nonumber  &&
H_n(\underline{g},\underline{f};\omega_1,\omega)
=\frac{1}{(n+1)!}\int_{u_j\in\textup{i}\R}
\prod_{1\leq i<j\leq n+1}{\gamma^{(2)}(g+u_i+u_j,f-u_i-u_j;\mathbf{\omega})\over
\gamma^{(2)}(u_i- u_j,u_j-u_i;\mathbf{\omega})}
\\ && \makebox[4em]{} \times
\prod_{j=1}^{n+1}\prod_{k=1}^{3}\gamma^{(2)}(g_k+ u_j,f_k-u_j;\mathbf{\omega})
\prod_{j=1}^n \frac{du_j}{\textup{i}\sqrt{\omega_1\omega_2}},
\end{eqnarray}
where $\sum_{j=1}^{n+1}u_j=0$,  $(n-1)(f+g)+\sum_{k=1}^3(g_k+f_k)=\omega_1+\omega_2$.
Then the hyperbolic reduction of the identity \eqref{AIIa-odd}  valid for odd $n$ has the form
\begin{eqnarray} \nonumber &&
H_n(\underline{g},\underline{f};\omega_1,\omega)
=\gamma^{(2)}\left(\tfrac{n+1}{2}g,\tfrac{n+1}{2}f;\mathbf{\omega}\right)
\prod_{1\leq i<k\leq 3}\gamma^{(2)}\left(\tfrac{n-1}{2}g+g_i+g_k,\tfrac{n-1}{2}f+f_i+f_k;\mathbf{\omega}\right)
\\ && \makebox[4em]{} \times
\prod_{j=1}^{(n+1)/2}\prod_{i,k=1}^3\gamma^{(2)}((j-1)(g+f)+g_i+f_k;\mathbf{\omega})
\\ &&\makebox[-2em]{}  \times
\prod_{j=1}^{(n-1)/2} \gamma^{(2)}(j(g+f);\mathbf{\omega})
\prod_{1\leq i<k\leq 3}\gamma^{(2)}((j-1)g+jf+g_i+g_k, jg+(j-1)f+f_i+f_k;\mathbf{\omega}).
\nonumber\end{eqnarray}
The hyperbolic reduction of identity \eqref{AIIa-even}  valid for even $n$ has the form
\begin{eqnarray} \nonumber && \makebox[-2em]{}
H_n(\underline{g},\underline{f};\omega_1,\omega)=\prod_{i=1}^{3}
\gamma^{(2)} (\tfrac{n}{2}g+g_i,\tfrac{n}{2}f+f_i;\mathbf{\omega})
\gamma^{(2)} ((\tfrac{n}{2}-1 )g+\sum_{k=1}^3 g_k,
(\tfrac{n}{2}-1 )f+\sum_{k=1}^3f_k;\mathbf{\omega})
 \\  &&\makebox[4em]{} \times
 \prod_{j=1}^{n/2}\Bigl(\gamma^{(2)}(j(g+f);\mathbf{\omega})
\prod_{i,k=1}^{3}\gamma^{(2)}((j-1)(g+f)+g_i+f_k);\mathbf{\omega})
\nonumber \\ &&  \makebox[2em]{} \times
\prod_{1\leq i<k\leq 3}\gamma^{(2)}((j-1)g+jf+g_i+g_k, jg+(j-1)f+f_i+f_k;\mathbf{\omega})\Bigr).
\end{eqnarray}

We use now the parametrization \eqref{zkomA} and \eqref{gkomA} and additionally set
\beq\label{gkomB}
g=\textup{i}\sqrt{\omega_1\omega_2}({N}+\delta {a}),\qquad
f=\textup{i}\sqrt{\omega_1\omega_2}({M}+\delta {b}),
\ee
where $M, N\in\Z$ and $a, b\in\CC$.
Then in the limit $\delta\to 0^+$, we come to  the function
\begin{eqnarray} \nonumber &&
R_n(\underline{a},\underline{N},\underline{b},\underline{M})
= \frac{1}{(2^{n+3}\pi)^n(n+1)!} \sum_{m_j\in \Z+\nu}\int_{x_j\in\textup{i}\R}
\prod_{1\leq j<k\leq n+1}\left[(x_j-x_k)^2+(m_j-m_k)^2\right]
\\ \nonumber && \makebox[4em]{} \times
\prod_{1\leq i<j\leq n+1}{\bf \Gamma}({a}+x_i+ x_j,{N}+m_i+ m_j){\bf \Gamma}({b}-x_i- x_j,{M}-m_i- m_j)
\\  && \makebox[4em]{} \times
\prod_{j=1}^{n+1}\prod_{k=1}^{3}{\bf \Gamma}(a_k+ x_j,N_k+ m_j){\bf \Gamma}(b_k- x_j,M_k- m_j)
\prod_{j=1}^n dx_j,
\label{rnn}\end{eqnarray}
where $\sum_{j=1}^{n+1}x_j=0$, $\sum_{j=1}^{n+1}m_j=0$ (for even $n$ this requires $\nu=0$),
and the balancing condition has the form
$$
(n-1)({N}+{M})+\sum_{k=1}^3(N_k+M_k)=0, \qquad (n-1)({a}+{b})+\sum_{k=1}^3(a_k+b_k)=-2\textup{i}.
$$
Convergence condition of the sum of integrals \eqref{rnn} is similar to the previous cases.
For large $|z_j|=\tfrac12|m_j+\textup{i}x_j|$ with fixed $j$, we have the asymptotics of
the integrand modulus
$$
\propto |z_j|^{2\textup{i}(a+b)(n-1)+2\sum_{k=1}^3 (\textup{i}(a_k+b_k)-2)+2}=|z_j|^{-6}
$$
due to the balancing condition, i.e. the convergence is guaranteed.

Finally, for odd $n$ we obtain  the identity
\begin{eqnarray} \nonumber &&
R_n(\underline{a},\underline{N},\underline{b},\underline{M})
= e^{\frac{\pi\textup{i}}{2}\varphi_{odd}}\,
{\bf \Gamma}\left(\tfrac{n+1}{2}{a}, \tfrac{n+1}{2}{N}\right){\bf \Gamma}
(\tfrac{n+1}{2}{b}, \tfrac{n+1}{2}{M})
  \\ \nonumber &&  \makebox[-2em]{} \times
\prod_{1\leq i<k\leq 3}{\bf \Gamma}\left(\tfrac{n-1}{2}{a}+a_i+a_k,\tfrac{n-1}{2}{N}+N_i+N_k\right)
{\bf \Gamma}\left(\tfrac{n-1}{2}{b}+b_i+b_k,\tfrac{n-1}{2}{M}+M_i+M_k\right)
\\ \label{phi_odd}  &&  \times
\prod_{j=1}^{(n+1)/2}\prod_{i,k=1}^3{\bf \Gamma}((j-1)({a}+{b})+a_i+b_k,(j-1)({N}+{M})+N_i+M_k)
\\ \nonumber && \times
\prod_{j=1}^{(n-1)/2}\Bigl(
\prod_{1\leq i<k\leq 3}\big[{\bf \Gamma}((j-1){a}+j{b}+a_i+a_k, (j-1){N}+j{M}+N_i+N_k)
\\  &&  \makebox[-2em]{} \times
 {\bf \Gamma}(j{a}+(j-1){b}+b_i+b_k,j{N}+(j-1){M}+M_i+M_k)\big]{\bf\Gamma}(j({a}+{b}),j({N}+{M}))\Bigr).
 \nonumber\end{eqnarray}
For even $n$ we have
\begin{eqnarray} \nonumber &&
R_n(\underline{a},\underline{N},\underline{b},\underline{M})
=e^{\frac{\pi\textup{i}}{2}\varphi_{even}}\, \prod_{i=1}^{3}{\bf \Gamma}\left(\tfrac{n}{2}{a}+a_i, \tfrac{n}{2}{N}+N_i\right){\bf \Gamma}\left(\tfrac{n}{2}{b}+b_i, \tfrac{n}{2}{M}+M_i\right)
 \\ \nonumber &&  \makebox[-2em]{} \times
 {\bf \Gamma}((\tfrac{n}{2}-1){a}+\sum_{k=1}^3a_k,(\tfrac{n}{2}-1){N}+\sum_{k=1}^3N_k)
 {\bf \Gamma}((\tfrac{n}{2}-1){b}+\sum_{k=1}^3b_k,(\tfrac{n}{2}-1){M}+\sum_{k=1}^3M_k)
\\  \label{phi_even}  &&  \times
\prod_{j=1}^{n/2}\Bigl(
\prod_{i,k=1}^{3} {\bf \Gamma}((j-1)({a}+{b})+a_i+b_k,(j-1)({N}+{M}))+N_i+M_k)
\\ \nonumber &&   \makebox[2em]{} \times
\prod_{1\leq i<k\leq 3}\big[{\bf \Gamma}((j-1){a}+j{b}+a_i+a_k, (j-1){N}+j{M}+N_i+N_k)
\\ &&  \makebox[-2em]{} \times
{\bf \Gamma}(j{a}+(j-1){b}+b_i+b_k,j{N}+(j-1){M}+M_i+M_k)\big]
{\bf \Gamma}(j({a}+{b}),j({N}+{M}))\Bigr).
 \nonumber\end{eqnarray}

Initial explicit expressions for the phases $\varphi_{odd}$ and $\varphi_{even}$
involve complicated quadratic polynomials of the discrete variables $N, N_k, M, M_k$.
However, they can be substantially simplified in the odd $n$ case to
\begin{equation}
\varphi_{\rm odd}=\frac{n}{4} (n^2-1)(N+M).
\label{phase_odd}\end{equation}
For even $n$ one has $\nu=0$ and we obtain the expression
\bea &&
\varphi_{\rm even}=(N^2+M^2)\left(-\tfrac14 n^2+\tfrac12 n+1\right)-2N\sum_{i=1}^3N_i-2M\sum_{i=1}^3M_i
\nonumber \\ &&   \makebox[3em]{} +\left(\sum_{i=1}^3N_i\right)^2+\left(\sum_{i=1}^3M_i\right)^2+n\sum_{i=1}^3(N_i^2+M_i^2).
\label{phase_even}\eea
One can check that both $\varphi_{\rm odd}$ and $\varphi_{\rm even}$ are even integer numbers.

In the following we consider limiting forms of some integrals described until now, excluding the type II
integrals presented in this section.

\section{Limiting cases of the beta integrals}

In this section we describe three more degenerations of the described beta integrals.
First we take relation (\ref{hyperAI}), shift parameters $g_l\to g_l-L,\, f_l\to f_l+L$
for all $l=1,\ldots, n+2$, and take the limit $L\to \textup{i}\infty$.
The resulting identity has the following form
\begin{eqnarray}
\frac{1}{n!}\int_{u_j\in\textup{i}\R }
\frac{\prod_{j=1}^{n}\prod_{\ell=1}^{n+2}\gamma^{(2)}(g_\ell+ u_j,f_\ell-u_j;\mathbf{\omega}) }
{\prod_{1\leq j<k\leq n}\gamma^{(2)}(\pm(u_j- u_k);\mathbf{\omega})
 } \prod_{j=1}^n \frac{du_j}{\textup{i}\sqrt{\omega_1\omega_2}}
= \prod_{\ell,s=1}^{n+2}\gamma^{(2)}(g_\ell+f_s;\mathbf{\omega})
\label{hyperAI2}\end{eqnarray}
with the same balancing condition $\sum_{l=1}^{n+2}(g_{\ell}+f_\ell)=Q$.

To prove this formula one should use the asymptotics \eqref{asy1}.
Under the taken conditions for parameters,
main contribution to the integral (\ref{hyperAI}) comes from the regions of large
$|u_j|,\, j=1,\ldots, n+1,$ which compensate the growth of $g_l, f_l$.
These regions are reached by the shifts $u_k\to u_k-nL,$ for some fixed $k$,
and $u_j\to u_j+L$ for all other variables, $j\neq k$. Because of the permutation symmetry in $u_j$,
all these $n+1$ arising asymptotic domains of integration are equivalent.
Therefore we consider only the integral emerging from the choice
$u_j\to u_j+L,\, j=1,\ldots, n$, $u_{n+1}\to u_{n+1}-nL$.
Then the terms whose $\gamma^{(2)}$-function arguments do not depend on $L$ remain intact
and one can see that they are included in (\ref{hyperAI2}).
The terms depending on $L$ will be replaced by the exponents of $B_{2,2}$-polynomials.
There is equal number of $+B_{2,2}$ and $-B_{2,2}$ terms in the exponents of asymptotics,
i.e. we can drop their constant terms.
Omitting the  factor ${\textup{i} \pi\over 2\omega_1\omega_2}$ common for all terms,
one can write for these exponents
\bea\nonumber &&
\sum_{l=1}^{n+2}\Big[\left(g_l+u_{n+1}-(n+1)L-{Q\over 2}\right)^2-\left(f_l-u_{n+1}+(n+1)L-{Q\over 2}\right)^2\Big]
\\ \label{hp1} &&  \makebox[1em]{}
+\sum_{j=1}^{n}\Big[\left(u_j-u_{n+1}+(n+1)L-{Q\over 2}\right)^2-\left(-u_j+u_{n+1}-(n+1)L-{Q\over 2}\right)^2\Big].
\eea
Taking into account the balancing condition and equality $\sum_{l=1}^{n}u_j=-u_{n+1}$,
we simplify this expression to
\beq\label{lfsc}
\sum_{l=1}^{n+2}\left(g_l^2-f_l^2\right)-Q\sum_{l=1}^{n+2}\left(g_l-f_l\right)+2LQ(n+1).
\ee
On the right-hand side we have
\bea \nonumber &&
\sum_{l=1}^{n+2}\Big[\left(G-g_l-(n+1)L-{Q\over 2}\right)^2-\left(F-f_l+(n+1)L-{Q\over 2}\right)^2\Big]
\eea
Using again the balancing condition $G+F=Q$, we obtain exactly the same expression (\ref{lfsc}).
Thus, we have no contribution from the divergent terms. Finally, reminding that there are $n+1$
domains of large $u_j$-variables yielding the same result in the large $L$ limit, as described
above, we replace the combinatorial factor $1/(n+1)!$ by $1/n!$.

Parametrizing now $g_l, f_l,$ and $u_j$ in (\ref{hyperAI2}) as in \eqref{zkomA}  and \eqref{gkomA},
we compute the limit $\delta\to 0^+$ and derive the identity
\begin{eqnarray} \nonumber &&  \makebox[-2em]{}
\frac{1}{(2^{n+1}\pi)^n n!}\sum_{m_j\in \Z+\nu}\int_{x_j\in \R }e^{\textup{i} \pi(n+2)\sum_{j=1}^n m_j}
\prod_{j=1}^{n} \prod_{\ell=1}^{n+2}{\bf \Gamma}(a_\ell+ x_j,N_\ell+ m_j){\bf \Gamma}(b_\ell- x_j,M_\ell- m_j)
\\ \nonumber && \makebox[4em]{}  \times
\prod_{1\leq j<k\leq n}\left[(x_j-x_k)^2+(m_j-m_k)^2\right] \prod_{j=1}^n dx_j
\\  &&  \makebox[6em]{}
=e^{\textup{i} \pi M}e^{\textup{i} \pi ((n+1)n-2)\nu}
\prod_{\ell,s=1}^{n+2}{\bf \Gamma}(a_\ell+b_s,N_\ell+M_s),
\label{complexAI1}\end{eqnarray}
where $\sum_{l=1}^{n+2}(N_l+M_l)=0$, $\sum_{l=1}^{n+2}(a_l+b_l)=-2\textup{i}$,
and  $M=\sum_{l=1}^{n+2}M_l$.
For $n=1$ this formula coincides with \eqref{MB}.
This formula was proved in \cite{DM2019} by a brute force method (see formula (3.5) there).
The parameter $\nu$ is redundant, since for $\nu=1/2$ the shifts $m_j\to m_j+1/2,\,
N_k\to N_k-1/2,\, M_k \to M_k+1/2$ yield the formula equivalent to the $\nu=0$ case.
As to the convergence condition in \eqref{complexAI1}, it is similar to the case \eqref{MB} since
the integrand modulus asymptotics
$$
\propto |z_j|^{\sum_{\ell=1}^{n+2} (\textup{i} (a_\ell+b_\ell)-2)+2(n-1)}=|z_j|^{-4},
\quad z_j=\tfrac12(m_j+\textup{i}x_j)\to\infty,
$$
guarantees validity of the relation

As a second example, we consider the reduction of relation  (\ref{hyperAI2}) following
from taking $g_{n+2}\to \textup{i}\infty$ jointly with the limit $f_{n+2}\to -\textup{i}\infty$
emerging through the balancing condition. In this way, one obtains the following exact formula
\begin{eqnarray} \nonumber &&
\frac{1}{n!}\int_{u_j\in\textup{i}\R }e^{{\textup{i} \pi\over \omega_1\omega_2}\left(\sum_{s=1}^{n+1}(g_s+f_s)\right)\left(\sum_{j=1}^{n}u_j\right)}
\,\frac{\prod_{j=1}^{n}\prod_{\ell=1}^{n+1}\gamma^{(2)}(g_\ell+ u_j,f_\ell-u_j;\mathbf{\omega}) }
{\prod_{1\leq j<k\leq n}\gamma^{(2)}(\pm(u_j- u_k);\mathbf{\omega})
 } \prod_{j=1}^n \frac{du_j}{\textup{i}\sqrt{\omega_1\omega_2}}
\\  && \makebox[2em]{}
=e^{{\textup{i} \pi\over 2\omega_1\omega_2}\Big[\left(\sum_{s=1}^{n+1}(g_s+f_s)\right)\left(\sum_{s=1}^{n+1}(f_s-g_s)\right)+\sum_{s=1}^{n+1}(g_s^2-f_s^2)\Big]}
{\prod_{\ell,s=1}^{n+1}\gamma^{(2)}(g_\ell+f_s;\mathbf{\omega})\over \gamma^{(2)}(\sum_{s=1}^{n+1}(g_s+f_s);\mathbf{\omega})}.
\label{hyperAIK}\end{eqnarray}
Using again the parametrization \eqref{zkomA} and   \eqref{gkomA}, we come in the limit
$\delta\to 0^+$ to the identity
\begin{eqnarray} \nonumber &&
\frac{1}{(2^{n+1}\pi)^n n!}\sum_{m_j\in \Z+\nu}\int_{x_j\in \R }e^{\textup{i} \pi(n+1)\sum_{j=1}^n m_j}
\prod_{j=1}^{n} \prod_{\ell=1}^{n+1}\Big[ {\bf \Gamma}(a_\ell+ x_j,N_\ell+ m_j)
\\ \nonumber && \makebox[4em]{}  \times
{\bf \Gamma}(b_\ell- x_j,M_\ell- m_j) \Big]
\prod_{1\leq j<k\leq n}\left[(x_j-x_k)^2+(m_j-m_k)^2\right] \prod_{j=1}^n dx_j
\\  && \makebox[2em]{}
=e^{\textup{i} \pi n(n+1)\nu}
{\prod_{\ell,s=1}^{n+1}{\bf \Gamma}(a_\ell+b_s,N_\ell+M_s)\over {\bf \Gamma}(\sum_{s=1}^{n+1}(a_s+b_s),\sum_{s=1}^{n+1}(N_s+M_s))}.
\label{complexAI3}\end{eqnarray}
This formula was obtained earlier by a different method in \cite{DM2019} (see formula (2.3a) there).
It is easy to see that for $\nu=1/2$ this formula becomes equivalent to the $\nu=0$ case
after the shifts $m_j\to m_j+1/2,\, N_k\to N_k-1/2$, and $M_k\to M_k+1/2$, i.e. we can set $\nu=0$.

In difference from the cases considered until now, the convergence condition of the
sum of integrals in  \eqref{complexAI3} imposes an additional constraint on the parameter values.
Namely, we have the integrand modulus asymptotics
$\propto |z_j|^{-\sum_{\ell=1}^{n+1} \textup{Im} (a_\ell+b_\ell)-4}.$
Therefore we have to demand that $\sum_{\ell=1}^{n+1} \textup{Im}(a_\ell+b_\ell)> -2.$

The third degeneration example deals with the identity (\ref{hyperCI}), where we take the limit
$g_{2n+3}\to +\textup{i}\infty$ in such a way that simultaneously $g_{2n+4}\to -\textup{i}\infty$ through the balancing condition.
This reduces equality  (\ref{hyperCI})  to the relation
\begin{eqnarray} \nonumber &&
\frac{1}{n!}\int_{u_j\in\textup{i}\R }\prod_{1\leq j<k\leq n}\frac{1}{\gamma^{(2)}(\pm u_j\pm u_k;\mathbf{\omega})}
\prod_{j=1}^n\frac{\prod_{\ell=1}^{2n+2}\gamma^{(2)}(g_\ell\pm u_j;\mathbf{\omega}) }
{\gamma^{(2)}(\pm 2 u_j;\mathbf{\omega}) } \frac{du_j}{2\textup{i}\sqrt{\omega_1\omega_2}}
\\  && \makebox[4em]{}
={\prod_{1\leq \ell<s\leq 2n+2}\gamma^{(2)}(g_\ell+g_s;\mathbf{\omega})\over \gamma^{(2)}\left(\sum_{k+1}^{2n+2}g_k ;\mathbf{\omega}\right)},
\label{hyperCIF}\end{eqnarray}
which is a hyperbolic analogue of the multiple Askey--Wilson integral considered by
Gustafson \cite{gus:some1,gus:some2}.

Using the parametrizations (\ref{zkom}) and (\ref{gkom}), in the limit $\delta\to 0^+$ we obtain
\begin{eqnarray} \nonumber &&
\frac{1}{(2^{2n+1}\pi)^n n!}\sum_{m_j\in \Z+\nu}\int_{x_j\in \R }
\prod_{1\leq j<k\leq n}\left[(x_j\pm x_k)^2+(m_j\pm m_k)^2\right]
\\ \nonumber   && \makebox[4em]{} \times
\prod_{j=1}^n\Big[ (x_j^2+m_j^2)\prod_{\ell=1}^{2n+2}{\bf \Gamma}(a_\ell\pm x_j,N_\ell\pm m_j)\Big]dx_j
\\  && \makebox[2em]{}
=e^{-\textup{i} \pi n(n+1)\nu} \,\frac{\prod_{1\leq \ell<s\leq 2n+2}{\bf \Gamma}(a_\ell+a_s,N_\ell+N_s)}
{{\bf \Gamma}(\sum_{s=1}^{2n+2}a_s,\sum_{s=1}^{2n+2}N_s)},
\label{complexCI'}\end{eqnarray}
where the chosen contours of integration are valid for Im$(a_\ell)<0$.
This formula was also obtained in \cite{DM2019} by a different method (see formula (2.3b) there).

Similar to the previous case, convergence of the sum of integrals in \eqref{complexCI'}
imposes an additional constraint
on the parameter values. If $z_j=\tfrac12(m_j+\textup{i}x_j) \to \infty$ for a fixed $j$, then
the integrand modulus asymptotics
$\propto |z_j|^{-\sum_{\ell=1}^{2n+2} 2\textup{Im}( a_\ell)-6}$
leads to the restriction
$\sum_{\ell=1}^{2n+2} \textup{Im}( a_\ell)  >-2.$

\section{Transformation rules}

The higher order elliptic hypergeometric integrals having more parameters than the
elliptic beta integrals do not admit exact evaluations. However, they admit highly nontrivial
symmetry transformations the univariate prototype of which has been discovered in
\cite{spi:theta2}. For a description of the general type I integrals on the root system $C_n$
we take $z_1,\ldots,z_n$ lying on the unit circle $\T$ and parameters $t_1,\ldots,t_{2n+2m+4}, p, q\in\mathbb{C}^\times$,
satisfying restrictions $|p|, |q|, |t_j|<1$ and the balancing condition $\prod_{j=1}^{2n+2m+4}t_j=(pq)^{m+1}$.
Then the function of interest has the form
\begin{eqnarray}
&& I_n^{(m)}(\underline{t};C)=\kappa_n^C\int_{\T^n}\prod_{1\leq j<k\leq n}\frac{1}{\Gamma(z_j^{\pm 1} z_k^{\pm 1};p,q)}
\prod_{j=1}^n\frac{\prod_{\ell=1}^{2n+2m+4}\Gamma(t_\ell z_j^{\pm 1};p,q)}
{\Gamma(z_j^{\pm2};p,q)}\, \frac{dz_j}{z_j}.
\label{integralsCI}\end{eqnarray}
For $n=1, m=0$ this is the computable elliptic beta integral \cite{spi:umn}.
For $n=m=1$ this is the $V$-function representing an elliptic analogue of the Euler--Gauss hypergeometric function \cite{spi:essays}. In \cite{rai:trans}, Rains proved the following transformation formula:
\begin{equation}
I_n^{(m)}(t_1,\ldots,t_{2n+2m+4};C)=\prod_{1\leq r<s\leq 2n+2m+4}\Gamma(t_rt_s;p,q)\;
I_m^{(n)}\left(\frac{\sqrt{pq}}{t_1},\ldots,\frac{\sqrt{pq}}{t_{2n+2m+4}};C\right),
\label{trafo}\end{equation}
where one has on the right-hand side the same contours of integration $\T$ provided $|\sqrt{pq}/t_j|<1$.
It represents a direct generalization of a particular $W(E_7)$-group symmetry transformation for
the $V$-function. Note that relation  \eqref{trafo} connects integrals on roots systems of
different rank.

One can reduce identity \eqref{trafo} to the hyperbolic level \cite{rai:limits}:
\beq\label{reflt}
H_n^{(m)}(\underline{g})=\prod_{1\leq r< s \leq 2n+2m+4}\gamma^{(2)}(g_r+g_s;\omega_1,\omega_2)
H_m^{(n)}\left(\underline{\lambda}\right)\, ,\quad \lambda_{\ell}={\omega_1+\omega_2\over 2}-g_{\ell}\, ,
\ee
where we have the function
\begin{eqnarray}
H_n^{(m)}(\underline{g})=\frac{1}{n!}\int_{u_j\in\textup{i}\R}\prod_{1\leq j<k\leq n}
\frac{1}{\gamma^{(2)}(\pm u_j\pm u_k;\mathbf{\omega})}
\prod_{j=1}^n\frac{\prod_{\ell=1}^{2n+2m+4}\gamma^{(2)}(g_\ell\pm u_j;\mathbf{\omega}) }
{\gamma^{(2)}(\pm 2 u_j;\mathbf{\omega}) } \frac{du_j}{2\textup{i}\sqrt{\omega_1\omega_2}}
\end{eqnarray}
with the balancing condition:
\beq\label{balcon2}
\sum_{\ell=1}^{2m+2n+4} g_\ell=(m+1)(\omega_1+\omega_2).
\ee

Applying again the parametrization \eqref{zkom} and \eqref{gkom} we take the limit $\delta\to 0^+$.
Using the relation
$$
\tfrac{1}{2}(\omega_1+\omega_2)-g_{\ell}
=\textup{i}\sqrt{\omega_1\omega_2}(-N_{\ell}+\delta(-a_{\ell}-\textup{i}))+O(\delta^2),\quad \ell=1,\ldots,2n+2m+4,
$$
we reduce equality \eqref{trafo} to the relation (after dropping a common diverging factor)
\begin{equation}\label{RR}
R_n^{(m)}(\underline{a},\underline{N})=(-1)^{(n+m-1)(n+m-2)\nu}\prod_{1\leq \ell<s\leq 2n+2m+4}{\bf \Gamma}(a_\ell+a_s,N_\ell+N_s)R_m^{(n)}(\underline{\hat{a}},\underline{\hat{N}}),
\end{equation}
where we have the following complex hypergeometric integral on the root system $C_n$
\begin{eqnarray} \nonumber &&
R_n^{(m)}(\underline{a},\underline{N})=\frac{1}{(2^{2n+1}\pi)^n n!}\sum_{m_j\in \Z+\nu}\int_{x_j\in\R}
\prod_{1\leq j<k\leq n}\left[(x_j\pm x_k)^2+(m_j\pm m_k)^2\right]
\\ && \makebox[4em]{} \times
\prod_{j=1}^n (x_j^2+m_j^2)
\prod_{\ell=1}^{2n+2m+4}{\bf \Gamma}(a_\ell\pm x_j,N_\ell\pm m_j)dx_j
\label{complexCI22}\end{eqnarray}
and $\hat{a}_\ell=-\textup{i}-a_\ell,\, \hat{N}_\ell=-N_\ell.$
The balancing condition \eqref{balcon2} reduces to two constraints
\beq\label{xkn2k}
\sum_{\ell=1}^{2n+2m+4}a_\ell=-2\textup{i}(m+1),\qquad \sum_{\ell=1}^{2n+2m+4} N_\ell=0.
\ee
For $z_j=\tfrac12(m_j+\textup{i}x_j)\to\infty$ the integrands on both sides of \eqref{RR}
have the asymptotics of modulus $\propto |z_j|^{-6},$
which guarantees convergence of the integrals and their sum without additional restrictions.

We consider now symmetry transformations for some other multidimensional integrals.
Define the following type II elliptic hypergeometric integral on the root system $C_n$
\begin{equation}
V(\underline{t};t;p,q)= \kappa_n^C\int_{\T^n}\prod_{1\le j<k\le n}\!
\frac{\Gamma(tz_j^{\pm 1}z_k^{\pm1};p,q)}
{\Gamma(z_j^{\pm1}z_k^{\pm1};p,q)}
\prod_{j=1}^n \frac{\prod_{k=1}^8\Gamma(t_kz_j^{\pm1};p,q)}
{\Gamma(z_j^{\pm 2};p,q)}\frac{dz_j}{z_j},
\label{V_Rains}\end{equation}
where $|t|, |t_j|<1$ and $t^{2n-2}\prod_{j=1}^8t_j=p^2q^2$ (the balancing condition).
For $n=1$ this is the elliptic analogue of the Euler--Gauss $_2F_1$ hypergeometric function
introduced in \cite{spi:theta2}.
The higher dimensional extension \eqref{V_Rains} was introduced by Rains in \cite{rai:trans}
where the $W(E_7)$ group of symmetry transformations for it was established.
As shown in \cite{spi:cs}, the function $V(t_1,\ldots,t_8;t;p,q)$
emerges in the van Diejen quantum many-body system \cite{vD0}
under certain restrictions on the parameters (the balancing condition)
as a normalization of a special eigenfunction of the Hamiltonian.

The key nontrivial $W(E_7)$-group symmetry transformation can be written in the form
\begin{equation}
V(\underline{t};t;p,q)=\prod_{l=0}^{n-1}\prod_{1\leq j< k\leq 4}\Gamma(t^lt_jt_k;p,q)
\prod_{5\leq j< k\leq 8}\Gamma(t^lt_jt_k;p,q)V(\underline{s};t;p,q),
\end{equation}
where $|t|, |t_j|, |s_j|<1$ and
\[
\left\{
\begin{array}{cl}
s_j =v t_j ,&   j=1,2,3,4,  \\
s_j = v^{-1} t_j, &    j=5,6,7,8,
\end{array}
\right.
\quad v=\sqrt{\frac{pqt^{1-n}}{t_1t_2t_3t_4}}
=\sqrt{\frac{t_5t_6t_7t_8}{pqt^{1-n}}}.
\]
In the limit \eqref{parlim2} it reduces to the relation
\begin{eqnarray} &&
H(\underline{g};g;\omega)=\prod_{l=0}^{n-1}\Big [
\prod_{1\leq j< k \leq 4}\gamma^{(2)}(lg+g_j+g_k;\mathbf{\omega})
\nonumber  \\ && \makebox[4em]{} \times
\prod_{5\leq j< k \leq 8}\gamma^{(2)}(lg+g_j+g_k;\mathbf{\omega})\Big]
\, H(\underline{h};g;\omega),
\end{eqnarray}
where
\begin{eqnarray}\nonumber &&
H(\underline{g};g;\omega)=
\int_{-i\infty}^{i\infty}\prod_{1\leq j < k \leq n}{\gamma^{(2)}(g\pm u_j\pm u_k;\mathbf{\omega})\over \gamma^{(2)}(\pm u_j\pm u_k;\mathbf{\omega})}
\prod_{j=1}^n{\prod_1^8\gamma^{(2)}(g_k\pm u_j;\omega_1,\omega_2)\over \gamma^{(2)}(\pm 2u_j;\mathbf{\omega})}{du_j\over 2i\sqrt{\omega_1\omega_2}},
\nonumber\end{eqnarray}
with the balancing condition $(2n-2)g+\sum_{j=1}^8 g_j=2Q$ and
\begin{equation}
h_j=g_j+\xi, \quad h_{j+4}=g_{j+4}-\xi, \quad j=1,2,3,4,\quad
\xi={1\over 2}(Q+(1-n)g-\sum_{j=1}^4 g_j).
\label{xi}\end{equation}

Recall again the parametrization \eqref{zkom} and \eqref{gkom} and set additionally
\beq\label{gkomg}
g=\textup{i}\sqrt{\omega_1\omega_2}(L+\delta \gamma),\quad \gamma\in\CC,
\quad L\in\Z\, .
\ee
In the limit $\delta\to 0^+$ the balancing condition converts to
\begin{equation}
(2n-2)L+\sum_{k=1}^8N_k=0\, ,\quad (2n-2)\gamma +\sum_{k=1}^8a_k=-4\textup{i},
\end{equation}
and we obtain the equality
\begin{eqnarray} \nonumber &&
R(\gamma,L,\underline{a},\underline{N},\nu)=e^{\textup{i} \pi n\left(\sum_{i=1}^4N_j\right)}\prod_{l=0}^{n-1}
\prod_{1\leq j< k \leq 4}{\bf \Gamma}(l\gamma+a_j+a_k,lL+N_j+N_k)
\\ &&  \makebox[2em]{} \times
\prod_{5\leq j< k \leq 8}{\bf \Gamma}(l\gamma+a_j+a_k,lL+N_j+N_k\,
)R(\gamma,L,\underline{b},\underline{M},\mu),
\label{Rtrafo}\end{eqnarray}
where
\begin{eqnarray} \nonumber &&
R(\gamma,L,\underline{a},\underline{N},\nu)=\sum_{m_j\in \Z+\nu}\int_{u_j\in \R}
\prod_{1\le j<k\le n}\! \frac{{\bf \Gamma}(\gamma \pm u_j\pm u_k, L\pm m_j\pm m_k)}
{{\bf \Gamma}(\pm u_j\pm u_k, \pm m_j\pm m_k)}
\\ && \makebox[4em]{}
\times \prod_{j=1}^n \prod_{\ell=1}^8{\bf \Gamma}(a_\ell\pm u_j, N_\ell\pm m_j)
(u_j^2+m_j^2)\, \prod_{j=1}^n du_j
\label{rlgann}\end{eqnarray}
and
\begin{eqnarray} \nonumber &&
b_j=a_j-\textup{i}+X,\;  M_j=N_j+K, \;
b_{j+4}=a_{j+4}+\textup{i}-X, \; M_{j+4}=N_j-K,\; j=1,2,3,4,
\end{eqnarray}
with $K$ and $X$ entering the parametrization of the $\xi$-variable in \eqref{xi}:
$$
{\xi\over \textup{i}\sqrt{\omega_1\omega_2}}=K-\delta\big(\textup{i}-X\big),\quad
X:=\tfrac{1}{2}\big((1-n)\gamma-\sum_{j=1}^4 a_j\big),
\quad K:=\tfrac{1}{2}\big((1-n)L-\sum_{j=1}^4 N_j\big).
$$
If ~$K$ is an integer, then in \eqref{Rtrafo} one has $\mu=\nu$. If~$K$ is a half-integer, then $\mu\neq \nu$.
The contours of integration of integrals on both sides of equality \eqref{Rtrafo} are   real lines,
which is a valid choice for Im$(a_\ell)$, Im$(\gamma)$, Im$(b_\ell) <0$.
In the limit $z_j=\tfrac12(m_j+\textup{i} x_j)\to \infty$ for a fixed $j$, we find the asymptotics
of the integrand modulus in \eqref{rlgann} $\propto |z_j|^{-6}$ due to the balancing condition,
which guarantees needed convergence.

General elliptic hypergeometric integrals of type I on the root
system $A_n$ are defined as
\begin{eqnarray} \nonumber && \makebox[-1em]{}
I_{n}^{(m)}(s_1,\ldots,s_{n+m+2};t_1,\ldots,t_{n+m+2};A)
\\ &&   \makebox[0em]{}
= \kappa_{n}^A\int_{\T^n} \prod_{1\leq j<k\leq n+1}
\frac{1}{\Gamma(z_jz_k^{-1},z_j^{-1}z_k;p,q)}\prod_{j=1}^{n+1}
\prod_{l=1}^{n+m+2}\Gamma(s_lz_j^{-1},t_lz_j;p,q)\; \prod_{k=1}^n\frac{dz_k}{z_k},
 \label{newint}\end{eqnarray}
where $|t_j|, |s_j|<1$,
$$
\prod_{j=1}^{n+1}z_j=1,\qquad \prod_{l=1}^{n+m+2}s_lt_l=(pq)^{m+1},
$$
and one sets $I_{0}^{(m)}=\prod_{l=1}^{m+2}\Gamma(s_l,t_l;p,q)$.
Consider the symmetry transformation for $I_{n}^{(m)}$-integrals found by Rains in \cite{rai:trans}.
We denote $T=\prod_{j=1}^{n+m+2}t_j,\, S=\prod_{j=1}^{n+m+2}s_j$,
so that $ST=(pq)^{m+1}$, and restrict all $|t_k|,$ $|s_k|,$ $|T^{\frac{1}{m+1}}/t_k|$,
$|S^{\frac{1}{m+1}}/s_k|<1$. Then the following symmetry transformation is true:
\begin{eqnarray}\nonumber
&&
I_{n}^{(m)}(t_1,\ldots,t_{n+m+2};s_1,\ldots,s_{n+m+2} ;A)
=\prod_{j,k=1}^{n+m+2}\Gamma(t_js_k;p,q) \\ && \makebox[1em]{} \times
I_{m}^{(n)}\left(\frac{T^{\frac{1}{m+1}}}{t_1},\ldots,
\frac{T^{\frac{1}{m+1}}}{t_{n+m+2}};
\frac{S^{\frac{1}{m+1}}}{s_1},\ldots,\frac{S^{\frac{1}{m+1}}}{s_{n+m+2}};
A\right).
\label{rains2}\end{eqnarray}
This relation generalizes a particular $W(E_7)$ symmetry transformation of the $V$-function
to multivariate integrals and connects integrals on root systems of different rank.

The hyperbolic analogue of integral \eqref{newint}
\begin{eqnarray}
H_{n}^{(m)}(\underline{g},\underline{f};A)=\frac{1}{(n+1)!}\int_{u_j\in\textup{i}\R}
\frac{\prod_{j=1}^{n+1}\prod_{\ell=1}^{n+m+2}\gamma^{(2)}(g_\ell+ u_j,f_\ell-u_j;\mathbf{\omega}) }
{\prod_{1\leq j<k\leq n+1}\gamma^{(2)}(\pm(u_j- u_k);\mathbf{\omega})
 } \prod_{j=1}^n \frac{du_j}{\textup{i}\sqrt{\omega_1\omega_2}},
\end{eqnarray}
where $\sum_{j=1}^{n+1}u_j=0$, $X=\sum_{\ell=1}^{n+m+2}g_\ell$ and $Y=\sum_{\ell=1}^{n+m+2}f_\ell$
with $X+Y=(m+1)(\omega_1+\omega_2)$,
satisfies the following reduction of identity \eqref{rains2} \cite{rai:limits}:
\begin{equation}\label{hypertran}
H_{n}^{(m)}(\underline{g},\underline{f};A)=\prod_{j,k=1}^{n+m+2}\gamma^{(2)}(g_j+f_k;\omega_1,\omega_2)
H_{m}^{(n)}(\underline{c},\underline{d};A),
\end{equation}
where
$$
c_\ell={X\over m+1}-g_\ell,\qquad      d_\ell={Y\over m+1}-f_\ell, \quad \ell=1,\ldots, n+m+2.
$$

In order to reduce this identity to the complex hypergeometric level we parametrize
the variables $g_\ell$ and $f_\ell$ in the same way as in \eqref{gkomA}. However, the integration
variables on two sides are parametrized differently. On the left-hand side we use again \eqref{zkomA},
but for the right-hand side we fix
\beq\label{zkomA'}
u_j= \textup{i}\sqrt{\omega_1\omega_2}\Big(-\frac{W}{m+1}+m_j+\delta x_j\Big),
\quad x_j\in\CC,\quad m_j\in\Z+\nu,\; \nu=0, \frac{1}{2}, \; j=1,\ldots,m+1,
\ee
where $W=\sum_{\ell=1}^{n+m+2}N_\ell$. Then, in the limit $\delta\to 0^+$, we obtain
from \eqref{hypertran} the relation
\begin{equation} \label{trafoAnI}
R_{n}^{(m)}(\underline{a},\underline{N},\underline{b},\underline{M},\nu)
=e^{\textup{i} \pi \varphi(\nu)} \prod_{j,k=1}^{n+m+2}\Gamma(a_j+b_k,N_j+M_k)
\, \tilde R_{m}^{(n)}(\underline{a},\underline{N},\underline{b},\underline{M},\nu),
\end{equation}
where on the left-hand side we have the function
\begin{eqnarray}  \nonumber &&  \makebox[-2em]{}
R_{n}^{(m)}(\underline{a},\underline{N},\underline{b},\underline{M},\nu)
=\frac{1}{(2^{n+3}\pi)^n (n+1)!}\sum_{m_j\in \Z+\nu}\int_{x_j\in\R}
\prod_{j=1}^{n+1} \prod_{\ell=1}^{n+m+2}\Big[{\bf \Gamma}(a_\ell+ x_j,N_\ell+ m_j)
\\ && \makebox[-1em]{} \times
{\bf \Gamma}(b_\ell- x_j,M_\ell- m_j)\Big]
 \prod_{1\leq j<k\leq n+1}\left[(x_j-x_k)^2+(m_j-m_k)^2\right] \prod_{j=1}^n dx_j
\label{rab}\end{eqnarray}
with $N_\ell,\, M_\ell\in\Z+\nu$ and the conditions $\sum_{j=1}^{n+1}x_j=\sum_{j=1}^{n+1}m_j=0,$
which assume that for even $n$ one has $\nu=0$. Additionally, we have the balancing conditions
$$
\sum_{\ell=1}^{n+m+2}(N_\ell+M_\ell)=0, \quad C+D=-2\textup{i}(m+1),\quad
C=\sum_{\ell=1}^{n+m+2}a_\ell,\quad D=\sum_{\ell=1}^{n+m+2}b_\ell.
$$

On the right-hand side we have the function
\begin{eqnarray}  \nonumber &&  \makebox[-2em]{}
\tilde R_{m}^{(n)}(\underline{a},\underline{N},\underline{b},\underline{M},\nu)
=\frac{1}{(2^{m+3}\pi)^m (m+1)!}\sum_{m_j\in \Z+\nu}\int_{x_j\in\R}
\prod_{j=1}^{m+1}\prod_{\ell=1}^{n+m+2}\Big[{\bf \Gamma}(\tilde a_\ell+ x_j,-N_\ell+ m_j)
\\ && \label{rab'} \makebox[-1em]{}
\times {\bf \Gamma}(\tilde b_\ell- x_j,-M_\ell- m_j)\Big]
\prod_{1\leq j<k\leq m+1}\left[(x_j-x_k)^2+(m_j-m_k)^2\right] \prod_{j=1}^m dx_j
\end{eqnarray}
with the conditions $\sum_{j=1}^{m+1}x_j=0$ and $\sum_{j=1}^{m+1}m_j=W$. Here we have
$$
\tilde{a}_\ell={C\over m+1}-a_\ell,\qquad \tilde{b}_\ell={D\over m+1}-b_\ell.
$$
The factor $e^{\textup{i}\pi\varphi(\nu)}$ present in \eqref{trafoAnI} has the phases
\begin{equation}\label{phases_trafo}
\varphi(0)= (n+m+1)W,\quad
\varphi(\tfrac12)=\tfrac14 (3m-n+2)(n+m+2)+W^2.
\end{equation}
We remind that for even $n$ we must have $\nu=0$ and note that for $\nu=1/2$ the
phase $\varphi(\tfrac12)$ is an integer number for odd $n$ and both even and odd $m$.

\section{Limiting cases of the transformation rules}

We describe now two further reductions of the described symmetry transformations.
First, we reduce \eqref{reflt} by taking the limit
\beq
g_{2n+2m+3}\to \textup{i}\infty,\quad g_{2n+2m+4}=(m+1)Q-\sum_{k=1}^{2n+2m+3}g_k\to -\textup{i}\infty.
\ee
In this way we obtain
\beq\label{refltred}
\tilde{H}_n^{(m)}(\underline{g})={\prod_{1\leq r< s \leq 2n+2m+2}\gamma^{(2)}(g_r+g_s;\omega_1,\omega_2)
\over \gamma^{(2)}(\sum_{i=1}^{2n+2m+2}g_i-mQ)}
\tilde{H}_m^{(n)}\left(\underline{\lambda}\right),\quad \lambda_{\ell}={\omega_1+\omega_2\over 2}-g_{\ell},
\ee
where
\begin{eqnarray}\makebox[-1em]{}
\tilde{H}_n^{(m)}(\underline{g})=\frac{1}{n!}\int_{u_j\in\textup{i}\R}\prod_{1\leq j<k\leq n}\frac{1}{\gamma^{(2)}(\pm u_j\pm u_k;\mathbf{\omega})}
\prod_{j=1}^n\frac{\prod_{\ell=1}^{2n+2m+2}\gamma^{(2)}(g_\ell\pm u_j;\mathbf{\omega}) }
{\gamma^{(2)}(\pm 2 u_j;\mathbf{\omega}) } \frac{du_j}{2\textup{i}\sqrt{\omega_1\omega_2}}.
\end{eqnarray}
Inserting the parametrization (\ref{zkom}) and (\ref{gkom}) in \eqref{refltred}, in the limit $\delta\to 0^+$
we derive
\begin{eqnarray} \nonumber &&
\frac{1}{2^{(2n+1)n}\pi^nn!}\sum_{m_j\in \Z+\nu}\int_{x_j\in\R}
\prod_{1\leq j<k\leq n}\left[(x_j\pm x_k)^2+(m_j\pm m_k)^2\right]
\prod_{j=1}^n (x_j^2+m_j^2)
\\  \nonumber && \makebox[4em]{} \times
\prod_{j=1}^n\prod_{\ell=1}^{2n+2m+2}{\bf \Gamma}(a_\ell\pm x_j,N_\ell\pm m_j)dx_j
\\  \label{dm67} && \makebox[2em]{}
=\frac{e^{\textup{i} \pi (n+m)(n+m+1)\nu}}{2^{(2m+1)m}\pi^mm!}{\prod\limits_{1\leq j< k \leq 2n+2m+2}{\bf \Gamma}(a_j+a_k, N_j+N_k)
\over \Gamma\Big(\sum\limits_{j=1}^{2n+2m+2} a_j+2m{\rm i},\sum\limits_{j=1}^{2n+2m+2} N_j\Big)}
\\ \nonumber && \makebox[4em]{} \times
\sum_{m_j\in \Z+\nu}\int_{x_j\in\R}
\prod_{1\leq j<k\leq m}\left[(x_j\pm x_k)^2+(m_j\pm m_k)^2\right]
\prod_{j=1}^m\Big[ (x_j^2+m_j^2)
\\  \nonumber && \makebox[4em]{} \times
\prod_{\ell=1}^{2n+2m+2}{\bf \Gamma}(-\textup{i}-a_\ell\pm x_j,-N_\ell\pm m_j)\Big]dx_j.
\end{eqnarray}
Consider convergence conditions for the sum of integrals on the left-hand side of this
equality. For $z_j=\tfrac12(m_j+\textup{i} x_j)\to\infty$ we find the integrand modulus asymptotics
$\propto |z_j|^{-\sum_{\ell=1}^{2n+2m+2} 2\,\textup{Im}(a_\ell)-4m-6}.$
Therefore convergence is guaranteed only under the additional condition
$\sum_{\ell=1}^{2n+2m+2} \textup{Im}(a_\ell)> -2m-2$. For the sum of integrals on the
right-hand side we obtain a different constraint following from the replacements
$a_\ell \to - \textup{i}-a_\ell$ and permutation of $n$ and $m$, namely,  $\sum_{\ell=1}^{2n+2m+2} \textup{Im}(a_\ell)< -2m$.

Formula  \eqref{dm67} was conjectured by Derkachov and Manashov in \cite{DM2019} (see formula (6.7) there).
Because of the uniformity of our limits, our derivation rigorously proves this conjecture under the
described conditions of the convergence of integrals.

Second, we reduce further the identity \eqref{reflt}. Denote $u_j, \, j=1,\ldots, n,$ and
$u_j',\, j=1,\ldots, m,$ integration variables in the integrals $H_n^{(m)}(\underline{g})$
and $H_m^{(n)}(\underline{\lambda})$, respectively.
Now make the replacements  $u_j\to u_j-L, \, u_j'\to u_j'-L,$ and
$g_l\to g_l+L, \, g_{l+n+m+2}= f_{l}-L,$ $l=1,\ldots, n+m+2$,
after which we take the limit $L\to -\textup{i}\infty$.
Then the balancing condition \eqref{balcon2} takes the form
\beq\label{balconn}
\sum_{\ell=1}^{n+m+2} (g_\ell+ f_\ell)=(m+1)Q.
\ee
As a result, we obtain
\beq\label{refltt}
\hat{H}_n^{(m)}(\underline{h},\underline{f})=\prod_{r,s=1}^{n+m+2}\gamma^{(2)}(g_r+f_s;\omega_1,\omega_2)
\hat{H}_m^{(n)}\left(\underline{\gamma},\underline{\lambda}\right),
\ee
where $ \lambda_{\ell}={Q\over 2}-g_{\ell},\, \gamma_{\ell}={Q\over 2}-f_{\ell},$ and
\begin{eqnarray}
\hat{H}_n^{(m)}(\underline{g},\underline{f})=\frac{1}{n!}\int_{u_j\in\textup{i}\R }
\frac{\prod_{j=1}^{n}\prod_{\ell=1}^{n+m+2}\gamma^{(2)}(g_\ell+ u_j,f_\ell-u_j;\mathbf{\omega}) }
{\prod_{1\leq j<k\leq n}\gamma^{(2)}(\pm(u_j- u_k);\mathbf{\omega})
 } \prod_{j=1}^n \frac{du_j}{\textup{i}\sqrt{\omega_1\omega_2}}.
\end{eqnarray}

Substituting in \eqref{refltt} $f_{n+m+2}=(m+1)Q-\sum_{l=1}^{n+m+1}(g_l+f_l)-g_{n+m+2}$
and taking the limit $g_{n+m+2}\to \textup{i} \infty$, we obtain
\begin{eqnarray} \nonumber  &&
\frac{1}{n!}\int_{u_j\in\textup{i}\R }
\exp\Big[ {\textup{i} \pi\over \omega_1\omega_2}\big[mQ-\sum_{\ell=1}^{n+m+1}(g_\ell+f_\ell)\big]
 \sum_{j=1}^n u_j  \Big]
 \\ \nonumber && \makebox[4em]{} \times
\frac{\prod_{j=1}^{n}\prod_{\ell=1}^{n+m+1}\gamma^{(2)}(g_\ell+ u_j,f_\ell-u_j;\mathbf{\omega}) }
{\prod_{1\leq j<k\leq n}\gamma^{(2)}(\pm(u_j- u_k);\mathbf{\omega})}
\prod_{j=1}^n \frac{du_j}{\textup{i}\sqrt{\omega_1\omega_2}}
 \\ \nonumber && \makebox[2em]{}
 = \exp\Big[ {\textup{i} \pi\over 2\omega_1\omega_2}\big[mQ\left(\sum_{\ell=1}^{n+m+1}(f_\ell-g_\ell)\right)
+2\sum_{1\leq r< s \leq n+m+1}(g_rg_s-f_rf_s)\big]\Big]
 \\ \nonumber && \makebox[4em]{}  \times
{\prod_{r,s=1}^{n+m+1}\gamma^{(2)}(g_r+f_s;\omega_1,\omega_2)
\over \gamma^{(2)}(\sum_{\ell=1}^{n+m+1}(g_\ell+f_\ell)-mQ;\omega_1,\omega_2)}
\\   \nonumber &&   \makebox[4em]{}  \times
\frac{1}{m!}\int_{u_j\in\textup{i}\R }
\exp\Big[{\textup{i} \pi\over \omega_1\omega_2}\big[-(m+1)Q+\sum_{\ell=1}^{n+m+1}(g_\ell+f_\ell)\big]
\sum_{j=1}^m u_j\Big]
\\  &&   \makebox[6em]{}  \times
\frac{\prod_{j=1}^{m}\prod_{\ell=1}^{n+m+1}\gamma^{(2)}(\lambda_\ell- u_j,\gamma_\ell+u_j;\mathbf{\omega}) }
{\prod_{1\leq j<k\leq m}\gamma^{(2)}(\pm(u_j- u_k);\mathbf{\omega})
 } \prod_{j=1}^m \frac{du_j}{\textup{i}\sqrt{\omega_1\omega_2}}.
\label{ggg}\end{eqnarray}
Applying to this relation the parametrization \eqref{zkomA} and  \eqref{gkomA}, in the limit $\delta\to 0^+$
we obtain
\begin{eqnarray} \nonumber  &&
\frac{1}{2^{(n+1)n}\pi^n n!}\sum_{m_j\in \Z+\nu}\int_{x_j\in \R }e^{\textup{i} \pi(n+m+1)\sum_{j=1}^n m_j}
\prod_{1\leq j<k\leq n}\left[(x_j-x_k)^2+(m_j-m_k)^2\right]
\\ \nonumber && \makebox[4em]{}  \times
\prod_{j=1}^{n} \prod_{\ell=1}^{n+m+1}{\bf \Gamma}(a_\ell+ x_j,N_\ell+ m_j){\bf \Gamma}(b_\ell- x_j,M_\ell- m_j)
\prod_{j=1}^n dx_j
\\  \nonumber && \makebox[-1em]{}
=e^{\textup{i} \pi (n-m)(n+m+1)\nu}\frac{e^{\textup{i} \pi m\sum_{l=1}^{n+m+1}(N_l-M_l)}}{2^{(m+1)m}\pi^m m!}
{\prod_{\ell,s=1}^{n+m+1}{\bf \Gamma}(a_\ell+b_s,N_\ell+M_s)\over {\bf \Gamma}(\sum_{s=1}^{n+m+1}(a_s+b_s)+2im,\sum_{s=1}^{n+m+1}(N_s+M_s))}
\\ &&\times
\sum_{m_j\in \Z+\nu}\int_{x_j\in \R }e^{\textup{i} \pi(n+m+1)\sum_{j=1}^m m_j}
\prod_{1\leq j<k\leq m}\left[(x_j-x_k)^2+(m_j-m_k)^2\right]
\label{dm66} \\ && \makebox[2em]{}  \times
\prod_{j=1}^{m} \prod_{\ell=1}^{n+m+1}{\bf \Gamma}(-\textup{i}-a_\ell- x_j,-N_\ell- m_j)
{\bf \Gamma}(-\textup{i}-b_\ell+x_j,-M_\ell+ m_j) \prod_{j=1}^m dx_j.
 \nonumber\end{eqnarray}
The parameter $\nu$ becomes redundant: for $\nu=1/2$ this formula is equivalent to the $\nu=0$ case,
which can be seen after the shifts $m_j\to m_j+1/2$, $N_k \to N_k-1/2$, and ~$M_k\to M_k+1/2$.
Examining the convergence of the sum of integrals in this equality in the way similar to the
previous case, we find that it is guaranteed under the similar constraints on the parameters
$$
-2m-2< \sum_{\ell=1}^{n+m+1} \textup{Im}(a_\ell+b_\ell)<-2m.
$$
Formula \eqref{dm66} was conjectured by Derkachov and Manashov \cite{DM2019} (see formula (6.6) there)
and our degeneration limit (when the integrand has a uniform asymptotics on compacta with exponentially small
corrections) proves it under the taken constraints.

\section{Rational hypergeometric beta integrals}

Now we would like to describe a different limit for some of the above integrals
from the hyperbolic level to the rational one. Let us set
\begin{equation}\label{om1om22}
b=\sqrt{\omega_1\over \omega_2}=1+\textup{i}\delta, \quad \delta\to 0^+,
\end{equation}
which assumes that $\omega_1+\omega_2=2\sqrt{\omega_1\omega_2}+O(\delta^2)$.
In the conformal field theory framework, the limit $b\to 1$ corresponds to
the central charge $c\to 25$.
Again, for a special choice of the argument, the $\gamma^{(2)}$-function shows
a singular behaviour determined in \cite{GSVS1}
\begin{equation}
\gamma^{(2)}(\sqrt{\omega_1\omega_2}(n+y\delta);\omega_1,\omega_2)\underset{\delta\to 0^+}{=} {\rm e}^{-\frac{\pi \textup{i}}{2}(n-1)^2}
(4\pi\delta)^{n-1}\left(1-\frac{n+\textup{i}y}{2}\right)_{n-1},
\label{limit2}\end{equation}
where $n\in\Z$, $y\in\CC$, $(a)_n$ is the standard Pochhammer symbol: $(a)_0=1$ and
\[
(a)_n=\frac{\Gamma(a+n)}{\Gamma(a)}=
\begin{cases}
 a(a+1)\cdots(a+n-1), & \text{for} \ n>0,\quad \\
\dfrac{1}{(a-1)(a-2)\cdots(a+n)}, &\text{for} \ n<0.
\end{cases}
\]

In order to apply this limit to the hyperbolic integral (\ref{hyperCI}), we parametrize corresponding integration variables $u_j$ and the constants $g_j$ in the following way
$$
u_j=\sqrt{\omega_1\omega_2}(m_j+\textup{i}x_j\delta), \quad j=1,\ldots, n, \quad
g_\ell=\sqrt{\omega_1\omega_2}(N_\ell+\textup{i}a_\ell\delta), \quad \ell=1,\ldots, 2n+4,
$$
where $x_j, a_\ell\in \CC$ and $m_j, N_\ell\in\Z+\nu,\; \nu=0, \frac{1}{2}$.
The balancing condition \eqref{balcon}  in the limit $\delta\to 0^+$ yields the constraints
\beq\label{mu84'}
\sum_{\ell=1}^{2n+4} a_\ell=0, \qquad \sum_{\ell=1}^{2n+4} N_\ell=2.
\ee
Using the asymptotics \eqref{limit2}, we obtain
$$
\gamma^{(2)}(\pm 2u_j)\to  {1\over (4\pi\delta)^2}{1\over x_j^2-m_j^2},
$$
$$
\prod_{1\leq j<k\leq n}\gamma^{(2)}(\pm u_j \pm u_k)\to  {1\over (4\pi\delta)^{2n(n-1)}}
{2^{2n(n-1)}\over \prod_{1\leq j<k\leq n}[(x_j \pm x_k)^2-(m_j \pm m_k)^2]},
$$
$$
\prod_{j=1}^{n}\prod_{\ell=1}^{2n+4}\gamma^{(2)}(g_\ell\pm u_j;\mathbf{\omega})\to
\frac{e^{2\pi \textup{i} n\nu}}{(4\pi\delta)^{4n(n+1)}}
\prod_{j=1}^n\prod_{\ell=1}^{2n+4}\left(1+{a_\ell-N_\ell\pm (x_j- m_j)\over 2}\right)_{N_\ell\pm m_j-1},
$$
\begin{eqnarray*} &&
\prod_{1\leq \ell < s\leq 2n+4}
\gamma^{(2)}(g_\ell+g_s;\mathbf{\omega})\to \textup{i}^n (-1)^{(n+1)(n+2)\nu+n+1} (4\pi\delta)^{-n(2n+3)}
\\ && \makebox[8em]{} \times
\prod_{1\leq \ell < s\leq 2n+4}\left(1+{a_\ell+a_s-N_\ell-N_s\over 2}\right)_{N_\ell+N_s-1}.
\end{eqnarray*}

Inserting these expressions in  \eqref{hyperCI} and dropping some common diverging multiplier,
we obtain the following evaluation formula for a type I rational hypergeometric integral
on the root system $C_n$:
\begin{eqnarray} \nonumber &&
\frac{1}{(2^{2n+1}\pi\textup{i})^n n!}\sum_{m_j\in \Z+\nu}\int_{x_j\in\R}
\prod_{1\leq j<k\leq n}\left[(x_j\pm x_k)^2-(m_j\pm m_k)^2\right]
\prod_{j=1}^n (x_j^2-m_j^2)
\\   \nonumber && \makebox[8em]{} \times
\prod_{j=1}^n\prod_{\ell=1}^{2n+4}\left(1+{a_\ell-N_\ell\pm (x_j- m_j)\over 2}\right)_{N_\ell\pm m_j-1}dx_j
\\ &&
=(-1)^{(n-1)(n-2)\nu+n+1}\prod_{1\leq \ell<s\leq 2n+4}\left(1+{a_\ell+a_s-N_\ell-N_s\over 2}\right)_{N_\ell+N_s-1}.
\label{complexCII}\end{eqnarray}
Tracing locations of the original integration contours in  \eqref{hyperCI}, one can see that
the contours of integration in \eqref{complexCII}  in this limit should separate poles associated with the
Pochhammer symbols with indices $N_\ell+m_j-1$ and $N_\ell-m_j-1$ from each other. Clearly, for
sufficiently large values of $|m_j|$ all integrand singularities lie only from one side
of the integration contour for each integration variable, which means that the contribution of corresponding integrals
equals to zero. This means that for any choice of variables $a_\ell$ and $N_\ell$ there
are only finitely many non-zero terms in the series and all of them are equal to sums of
finitely many pole residues. In turn, this means that on the left hand-side of equality
\eqref{complexCII} one has a rational function of free parameters.

For demonstrating the nature of the derived identity, let us give two examples. First, we fix
$$
\nu=0, \qquad N_j=0, \quad j=1,2,\ldots, 2n+2,\quad N_{2n+3}=N_{2n+4}=1,
$$
and assume that Im$(a_j)<0,\, j=1,\ldots,2n+2$. Then, only the terms with $m_j=0$
can give nonzero contributions to the bilateral sums over $m_j$. Applying the residue
calculus to the emerging multiple integral of a rational function of $x_j$, we come
to the identity
\begin{eqnarray}\label{b1id1} &&
\frac{1}{(\pi\textup{i})^n n!}\int_{\R^n}
\frac{\prod_{1\leq j<k\leq n}(x_j^2- x_k^2)^2\prod_{j=1}^n x_j^2 dx_j}
{\prod_{j=1}^n\prod_{\ell=1}^{2n+2}(a_\ell^2-x_j^2)}
\\ && \makebox[1em]{}
=(-1)^{\frac{n(n-1)}{2}}\sum_{1\leq \ell_1<\ldots <\ell_n\leq 2n+2}
\frac{\prod_{j=1}^n a_{\ell_j}}{\prod_{j=1}^n\prod_{ \makebox[-2em]{} \ell=1 \atop \ell\neq \ell_1,\ldots, \ell_n}^{2n+2}
(a_\ell^2-a_{\ell_j}^2)}
 =\frac{(-1)^n\sum_{\ell=1}^{2n+2} a_\ell}{\prod_{1\leq \ell<k\leq 2n+2}(a_\ell+a_k)}
\nonumber\end{eqnarray}
after dropping the common multiplier $2^{2n^2+3n}$. For $n=1$ this identity was derived in \cite{GSVS2}.
The second expression is obtained by taking sequentially residues at the points $x_j=a_{\ell_k}$
after ordering the integration contours in one particular way.
This relation describes the leading $a_j\to 0$ asymptotics of the Gustafson
integral presented in Theorem 9.3 of \cite{gus:trams}.

In the second example, we fix
$$
\nu=0, \qquad N_j=0, \quad j=1,2,\ldots, 2n+3,\quad N_{2n+4}=2,
$$
and assume that Im$(a_j)<0,\, j=1,\ldots,2n+3$. Again, only the terms with $m_j=0$
can contribute to the sums over $m_j$. Applying the residue calculus  to the remaining
multiple integral of a rational function of $x_j$,  after dropping the common multiplier $2^{2n^2+3n}$
in the same way as in the previous example, we come to the identity
\begin{eqnarray}\label{b1id2} &&
\frac{1}{(\pi\textup{i})^n n!}\int_{\R^n}
\frac{\prod_{1\leq j<k\leq n}(x_j^2- x_k^2)^2\prod_{j=1}^n x_j^2 (a_{2n+4}^2-x_j^2)dx_j}
{\prod_{j=1}^n\prod_{\ell=1}^{2n+3}(a_\ell^2-x_j^2)}
\\ && \makebox[1em]{}
=(-1)^{\frac{n(n-1)}{2}}\sum_{1\leq \ell_1<\ldots <\ell_n\leq 2n+3}
\frac{\prod_{j=1}^n a_{\ell_j}(A^2-a_{\ell_j}^2)}
{\prod_{j=1}^n\prod_{ \makebox[-2em]{} \ell=1 \atop \ell\neq \ell_1,\ldots, \ell_n}^{2n+3}(a_\ell^2-a_{\ell_j}^2)}
 =\frac{(-1)^n\prod_{\ell=1}^{2n+3}(A- a_\ell)}{\prod_{1\leq \ell<k\leq 2n+3}(a_\ell+a_k)},
\nonumber\end{eqnarray}
where $A=\sum_{\ell=1}^{2n+3} a_\ell$. For $n=1$ this identity was also
derived in \cite{GSVS2}. This relation describes the leading $a_j\to 0$ asymptotics
of another Gustafson integral from \cite{gus:some2}, more complicated than in the previous case.
Note that formula \eqref{b1id2} reduces to \eqref{b1id1} in the limit $a_{2n+3}\to\infty$.

Now we turn to the integral \eqref{hyperAI}.
For $j=1,\ldots, n+1$ and $\ell=1,\ldots, 2n+4$ we parametrize
$$
u_j=\sqrt{\omega_1\omega_2}(m_j+\textup{i}x_j\delta), \;  \quad
g_\ell=\sqrt{\omega_1\omega_2}(N_\ell+\textup{i}a_\ell\delta), \;
f_\ell=\sqrt{\omega_1\omega_2}(M_\ell+\textup{i}b_\ell\delta),
$$
and take the limit $\delta\to 0^+$.
Here we have the variables $N_\ell, M_\ell, m_j \in\Z+\nu,\, \nu=0,\,\frac{1}{2}$, so that
$N_\ell\pm m_j$, $M_\ell\pm m_j$ and $N_\ell+M_s$, $N_\ell+N_s$, $M_\ell+M_s$  take integer values.
The balancing condition $G+F=\omega_1+\omega_2$  yields the constraints
\begin{eqnarray}  \nonumber &&
A+B=0,\quad A=\sum_{\ell=1}^{n+2}a_\ell,\quad B=\sum_{\ell=1}^{n+2}b_\ell,
\\&&
N+M=2,\quad N=\sum_{\ell=1}^{n+2} N_\ell,\quad M=\sum_{\ell=1}^{n+2} M_\ell.
\end{eqnarray}
Additionally, the condition $\sum_{j=1}^{n+1}u_j=0$ results in two more constraints
$\sum_{j=1}^{n+1}x_j=0,\, \sum_{j=1}^{n+1}m_j=0.$ The latter equality implies
that for even $n$ the discrete variables $m_j$ can take only integer values, i.e. the value $\nu=1/2$ is
forbidden. For odd $n$, both values $\nu=0,\tfrac12$ are allowed.

Using the asymptotics \eqref{limit2}, in the limit $\delta\to 0^+$ we obtain
\beq\label{limok}
\prod_{1\leq j<k\leq n+1}\gamma^{(2)}(\pm(u_j-u_k))\to
 {1\over (4\pi\delta)^{n(n+1)}}{2^{n(n+1)}\over \prod_{1\leq j<k\leq n+1}[(x_j-x_k)^2-(m_j-m_k)^2]},
\ee
\begin{eqnarray}\label{trfg}
&&\prod_{j=1}^{n+1}\prod_{\ell=1}^{n+2}\gamma^{(2)}(g_\ell+ u_j,f_\ell-u_j;\mathbf{\omega})\to
\frac{e^{-\pi\textup{i}(n+1)\left(\frac{1}{2}\sum_{\ell=1}^{n+2}(N_\ell^2+M_\ell^2)+(n+2)\nu^2+n\right)}}
{(4\pi\delta)^{2(n+1)^2}}\\ \nonumber
&&\times \prod_{j=1}^{n+1} \prod_{\ell=1}^{n+2}\left(1+{a_\ell+ x_j-(N_\ell+ m_j)\over 2}\right)_{N_\ell+ m_j-1}
 \left(1+{b_\ell- x_j-(M_\ell- m_j)\over 2}\right)_{M_\ell- m_j-1},
\end{eqnarray}
\begin{eqnarray}\label{tgyh}
&&\prod_{\ell=1}^{n+2}\gamma^{(2)}(G-g_\ell, F-f_\ell;\mathbf{\omega})
\prod_{\ell,s=1}^{n+2}\gamma^{(2)}(g_\ell+f_s;\mathbf{\omega})\to \\ \nonumber
&&{e^{-\textup{i}{\pi\over 2}\left( 2N^2(n-1)+n^2+2n+(n+3)\sum_{\ell=1}^{n+2}(N_\ell^2+M_\ell^2)\right)}\over (4\pi\delta)^{n^2+2n+2}}\prod_{\ell=1}^{n+2}
\left(1+{A-a_\ell-(N-N_\ell)\over 2}\right)_{N-N_\ell-1}\\ \nonumber
&&\times\prod_{\ell=1}^{n+2}\left(1+{B-b_\ell-(M-M_\ell)\over 2}\right)_{M-M_\ell-1}
\prod_{\ell,s=1}^{n+2}\left(1+{a_\ell+b_s-(N_\ell+M_s)\over 2}\right)_{N_\ell+M_s-1}.
\end{eqnarray}

Inserting \eqref{limok}, \eqref{trfg}, \eqref{tgyh} in equality \eqref{hyperAI} and dropping some common diverging multiplier, we obtain the following evaluation formula for a type I rational hypergeometric
integral on the root system $A_n$:
\begin{eqnarray} \nonumber &&
\frac{1}{(2^{n+3}\pi\textup{i})^n (n+1)!}\sum_{m_1,\ldots, m_n\in \Z+\nu}\int_{x_j\in\R}\prod_{1\leq j<k\leq n+1}\Big[(x_j-x_k)^2-(m_j-m_k)^2\Big]
\\ \nonumber && \makebox[-2em]{} \times
\prod_{j=1}^{n+1} \prod_{\ell=1}^{n+2}\Big(1+{a_\ell+ x_j-N_\ell- m_j\over 2}\Big)_{N_\ell+ m_j-1}
 \Big(1+{b_\ell- x_j-M_\ell+m_j\over 2}\Big)_{M_\ell- m_j-1}
\prod_{j=1}^n dx_j
\\  \label{AIIb1} && \makebox[2em]{}
=e^{-\textup{i} \pi(n-1)(\frac12 n+N^2-(n+2)\nu^2)}
\prod_{\ell=1}^{n+2}\Big(1+{A-a_\ell-N+N_\ell\over 2}\Big)_{N-N_\ell-1}
\\ && \times
\prod_{\ell=1}^{n+2}\Big(1+{B-b_\ell-M+M_\ell\over 2}\Big)_{M-M_\ell-1}
\prod_{\ell,s=1}^{n+2}\Big(1+{a_\ell+b_s-N_\ell-M_s\over 2}\Big)_{N_\ell+M_s-1}.
\nonumber\end{eqnarray}
In the same way as in the limit $b\to\textup{i}$, here we have the restrictions
Im$(a_\ell)$, Im$(b_\ell)<0$. The contours of integration (real lines) pass below
the poles stemming from the Pochhammer symbols $(\ldots)_{N_\ell+ m_j-1}$,
for  $j=1,\ldots, n$, and $(\ldots)_{M_\ell-m_{n+1}-1}$,
and above the poles coming from $(\ldots)_{M_\ell- m_j-1}$, for  $j=1,\ldots, n$,
and $(\ldots)_{N_\ell+m_{n+1}-1}$.

Take $m_j\to +\infty$ (or $-\infty$) for a fixed $j\neq n+1$, and assume that $m_{n+1}\to -\infty$
(or $+\infty$). Then, evidently, for sufficiently large $m_j$ the upper (or lower) half-plane for $x_j$
does not contain poles, i.e. the integral over $x_j$ vanishes. Suppose now that $m_j,\, m_{n+1}\to +\infty$.
This means that $m_i\to -\infty$ for some $i$, $n\geq i\neq j$. But then the lower half-plane for $x_i$
does not contain poles for sufficiently large $|m_i|$ and the integral over $x_i$ vanishes.
The third option is that $m_i$ compensates $m_j$ and $m_{n+1}$ is finite. In this case the
poles in the upper $x_j$ half-plane can emerge only from the Pochhammer symbol $(\ldots)_{M_\ell-m_{n+1}-1}$,
and the poles in the lower $x_i$ half-plane --- from  $(\ldots)_{N_\ell+m_{n+1}-1}$.
Taking residues of such poles in $x_j$-variable removes poles in the lower half-plane for the
$x_i$-variable and the corresponding integral vanishes. So, we have a finite sum of
integrals of rational functions representing a particular rational function of the external parameters.

Similar to the $C_n$-root system integrals, we give three simple examples of the realization of
\eqref{AIIb1}, all having $\nu=0$. In the first case set all $M_k=0$ and $N_\ell=0,\, \ell=1,\ldots, n,$
and $N_{n+1}=N_{n+2}=1$.
It is not difficult to see that for $m_j<0,\, m_{n+1}>0$ there are no
poles in the upper $x_j$-plane, and for $m_j,\, m_{n+1}<0$ there are no
poles in the upper $x_i$-plane for each $m_i>0$ (there is at least one such $m_i$
to compensate $m_j$). The same picture holds for $m_{n+1}=0$, unless all $m_j=0$. Therefore
only when all $m_j=0$ one gets a nonzero term in the sum which yields the following integral identity
\begin{eqnarray}\nonumber &&
\frac{1}{(2\pi\textup{i})^n (n+1)!}\int_{\R^n}
\frac{\prod_{1\leq j<k\leq n+1}(x_j- x_k)^2\, \prod_{j=1}^n dx_j}
{\prod_{j=1}^{n+1}\big[\prod_{\ell=1}^{n}(a_\ell+x_j)\prod_{r=1}^{n+2}(b_r-x_j)\big]}
\\ && \makebox[-1em]{}
=\sum_{s=0}^n\frac{s!(n-s)!}{(n+1)!}
\sum_{1\leq r_1<\ldots <r_s\leq n+2 \atop 1\leq \ell_1<\ldots <\ell_{n-s}\leq n}
\frac{(-1)^{\frac{s(s-1)}{2}+\frac{(n-s)(n-s-1)}{2}}}
{\prod_{j=1}^s\prod_{ \makebox[-2em]{} r=1 \atop r\neq r_1,\ldots, r_s}^{n+2}(b_r-b_{r_j})
\prod_{k=1}^{n-s}\prod_{ \makebox[-3em]{} \ell=1 \atop \ell\neq \ell_1,\ldots, \ell_{n-s}}^{n}(a_\ell-a_{\ell_k})}
\nonumber  \\  && \makebox[4em]{} \times
\frac{\prod_{j=1}^s (b_{r_j}+B_s-A_s)\prod_{k=1}^{n-s}(a_{\ell_k}+A_s-B_s)}
{\prod_{ \makebox[-2em]{} r=1 \atop r\neq r_1,\ldots, r_s}^{n+2} (b_r+B_s-A_s)
\prod_{ \makebox[-3em]{} \ell=1 \atop \ell \neq \ell_1,\ldots, \ell_{n-s}}^{n}(a_\ell +A_s-B_s)}
\nonumber  \\  && \makebox[4em]{} \times
\frac{\prod_{j=1}^{s}\prod_{k=1}^{n-s}(b_{r_j}+a_{\ell_k})^2}
{\prod_{j=1}^{s}\prod_{\ell=1}^{n}(b_{r_j}+a_{\ell})\prod_{k=1}^{n-s}\prod_{r=1}^{n+2}(b_r+a_{\ell_k})}
\nonumber  \\ && \makebox[6em]{}
 =\frac{(-1)^{\frac{n(n+1)}{2}} \prod_{\ell=1}^{n} (B+a_\ell)}
 {\prod_{k=1}^{n+2}(B-b_k)\prod_{\ell=1}^{n}\prod_{k=1}^{n+2} (a_\ell+b_k)}
\label{b1id3}\end{eqnarray}
after dropping the common multiplier $2^{n^2+2n+2}$. Here we have
$$
B=\sum_{\ell=1}^{n+2}b_\ell, \qquad A_s:= \sum_{k=1}^{n-s} a_{\ell_k}, \qquad  B_s:= \sum_{j=1}^s b_{r_j}.
$$

Next example uses the values $N_\ell=M_\ell=0,\, \ell=1,\ldots, n+1,$ and $N_{n+2}=M_{n+2}=1$. Again,
only the term with the summation indices $m_j=0$ gives a contribution yielding the identity
\begin{eqnarray}\nonumber &&
\frac{1}{(2\pi\textup{i})^n (n+1)!}\int_{\R^n}
\frac{\prod_{1\leq j<k\leq n+1}(x_j- x_k)^2\, \prod_{j=1}^n dx_j}
{\prod_{j=1}^{n+1}\big[\prod_{\ell=1}^{n+1}(a_\ell+x_j)\prod_{k=1}^{n+1}(b_k-x_j)\big]}
\\ && \makebox[-1em]{}
=\sum_{s=0}^n\frac{s!(n-s)!}{(n+1)!}
\sum_{1\leq r_1<\ldots <r_s\leq n+1 \atop 1\leq \ell_1<\ldots <\ell_{n-s}\leq n+1}
\frac{(-1)^{\frac{s(s-1)}{2}+\frac{(n-s)(n-s-1)}{2}}}
{\prod_{j=1}^s\prod_{ \makebox[-2em]{} r=1 \atop r\neq r_1,\ldots, r_s}^{n+1}(b_r-b_{r_j})
\prod_{k=1}^{n-s}\prod_{ \makebox[-3em]{} \ell=1 \atop \ell\neq \ell_1,\ldots, \ell_{n-s}}^{n+1}(a_\ell-a_{\ell_k})}
\nonumber  \\  && \makebox[4em]{} \times
\frac{\prod_{j=1}^s (b_{r_j}+B_s-A_s)\prod_{k=1}^{n-s}(a_{\ell_k}+A_s-B_s)}
{\prod_{ \makebox[-2em]{} r=1 \atop r\neq r_1,\ldots, r_s}^{n+1} (b_r+B_s-A_s)
\prod_{ \makebox[-3em]{} \ell=1 \atop \ell \neq \ell_1,\ldots, \ell_{n-s}}^{n+1}(a_\ell +A_s-B_s)}
\nonumber  \\  && \makebox[4em]{} \times
\frac{\prod_{j=1}^{s}\prod_{k=1}^{n-s}(b_{r_j}+a_{\ell_k})^2}
{\prod_{j=1}^{s}\prod_{\ell=1}^{n+1}(b_{r_j}+a_{\ell})\prod_{k=1}^{n-s}\prod_{r=1}^{n+1}(b_r+a_{\ell_k})}
\nonumber  \\ && \makebox[6em]{}
 =\frac{(-1)^{\frac{n(n+1)}{2}} \sum_{\ell=1}^{n+1} (a_\ell+b_\ell)}
 {(\sum_{\ell=1}^{n+1}a_\ell)(\sum_{k=1}^{n+1}b_k)\prod_{\ell,k=1}^{n+1}(a_\ell+b_k)}
\label{b1id4}\end{eqnarray}
after dropping the common multiplier $2^{n^2+2n+2}$. For $n=1$ both \eqref{b1id3} and \eqref{b1id4}
coincide with the identity \eqref{b1id1}.

Finally, for the choice $N_\ell=M_\ell=0, \ell=1,\ldots, n+1,$ and $N_{\ell+2}=2,\, M_{n+2}=0$
we find the relation
\begin{eqnarray}\nonumber &&
\frac{1}{(2\pi\textup{i})^n (n+1)!}\int_{\R^n}
\frac{\prod_{1\leq j<k\leq n+1}(x_j- x_k)^2\prod_{j=1}^{n+1}(-a_{n+2}-x_j)\, \prod_{j=1}^n dx_j}
{\prod_{j=1}^{n+1}\big[\prod_{\ell=1}^{n+1}(a_\ell+x_j)\prod_{k=1}^{n+2}(b_k-x_j)\big]}
\\ && \makebox[-2em]{}
=\sum_{s=0}^n\frac{s!(n-s)!}{(n+1)!}\makebox[-0.5em]{}
\sum_{1\leq r_1<\ldots <r_s\leq n+2 \atop 1\leq \ell_1<\ldots <\ell_{n-s}\leq n+1}
\frac{(-1)^{\frac{s(s-1)}{2}+\frac{(n-s)(n-s-1)}{2}}}
{\prod_{j=1}^s\prod_{ \makebox[-2em]{} r=1 \atop r\neq r_1,\ldots, r_s}^{n+2}(b_r-b_{r_j})
\prod_{k=1}^{n-s}\prod_{ \makebox[-3em]{} \ell=1 \atop \ell\neq \ell_1,\ldots, \ell_{n-s}}^{n+1}(a_\ell-a_{\ell_k})}
\nonumber  \\  && \makebox[4em]{} \times
\frac{(B+\sum_{\ell=1}^{n+1}a_\ell+B_s-A_s)\prod_{j=1}^s (b_{r_j}+B_s-A_s)\prod_{k=1}^{n-s}(a_{\ell_k}+A_s-B_s)}
{\prod_{ \makebox[-2em]{} r=1 \atop r\neq r_1,\ldots, r_s}^{n+2} (b_r+B_s-A_s)
\prod_{ \makebox[-3em]{} \ell=1 \atop \ell \neq \ell_1,\ldots, \ell_{n-s}}^{n+1}(a_\ell +A_s-B_s)}
\nonumber  \\  && \makebox[1em]{} \times
\frac{\prod_{j=1}^{s}\prod_{k=1}^{n-s}(b_{r_j}+a_{\ell_k})^2\prod_{j=1}^s(B+\sum_{\ell=1}^{n+1}a_\ell-b_{r_j})
\prod_{k=1}^{n-s}(B+\sum_{\ell=1}^{n+1}a_\ell+a_{\ell_k})}
{\prod_{j=1}^{s}\prod_{\ell=1}^{n+1}(b_{r_j}+a_{\ell})\prod_{k=1}^{n-s}\prod_{r=1}^{n+2}(b_r+a_{\ell_k})}
\nonumber  \\ && \makebox[6em]{}
 =\frac{(-1)^{\frac{n(n+1)}{2}} \prod_{\ell=1}^{n+1}(a_\ell+B)\prod_{k=1}^{n+2}(\sum_{i=1}^{n+1}a_i+B- b_k)}
 {(\sum_{\ell=1}^{n+1}a_\ell)\prod_{k=1}^{n+2}(B-b_k)\prod_{\ell=1}^{n+1}\prod_{k=1}^{n+2}(a_\ell+b_k)}
\label{b1id5}\end{eqnarray}
after dropping the common multiplier $2^{n^2+2n+2}$. For $n=1$ this identity coincides with \eqref{b1id2},
it reduces to \eqref{b1id3} in the limit $a_{n+1}\to\infty$
and to \eqref{b1id4} in the limit $b_{n+2}\to\infty$.

In a similar way, one can apply the described special $b\to1$ limit to all other hyperbolic
hypergeometric identities described in the present paper, as well as to those following
from a large number of conjectural relations for elliptic hypergeometric integrals
(superconformal indices of four-dimensional quantum field theory models) described in \cite{SV}.
This would lead to rational hypergeometric functions on root systems with nontrivial symmetry
properties. In the present paper we have chosen only the simplest examples of such functions
in order to illustrate the structure of emerging objects and identities.

\section{Conclusion}

In the present paper we have systematically investigated multivariate complex hypergeometric
functions on root systems. By taking appropriate  limits in the most popular elliptic hypergeometric integrals
on root systems $A_n$ and $C_n$, we arrived at bilateral infinite sums of
multiple integrals of the Mellin--Barnes type.  As to the identities for integrals over the complex planes
generalizing the complex beta-integral \eqref{STR}, they can be obtained by applying the quasiclassical
limit to obtained relations, as demonstrated in \cite{DM2019} on some particular examples.
It is necessary to investigate such a limit for all obtained identities with a view on
possible applications to two-dimensional conformal field theory. We described also two multiple rational
hypergeometric beta integrals on the root systems  $A_n$ and $C_n$  admitting exact evaluations and
gave a number of examples of corresponding rational function identities.

Special functions of hypergeometric type are known to be related to the representation
theory of Lie algebras or groups and their deformations. Considerations of the complex hypergeometric functions
are lying far behind the forefront of investigations in this direction. Basically, only $3j$- and
$6j$-symbols for the group $SL(2,\mathbb{C})$ have been related to such hypergeometric functions
\cite{DS2017,Ismag2,Naimark}. For the most complicated identities, there emerges the additional
discrete parameter $\nu=0, \tfrac12$ whose group-theoretical meaning is not clear
at the moment. In particular, it is not clear also what kind of algebraic structures are hidden
behind the $b\to 1$ limiting relations for the hyperbolic hypergeometric identities
described in the previous section.

Elliptic hypergeometric integrals unify all special functions of hypergeometric type, i.e.
they are universal special functions. Their most beautiful examples and identities describe superconformal
indices of four-dimensional supersymmetric field theories and corresponding Seiberg dualities \cite{DO,SV}.
Degeneration of these functions is related to dimensional reductions of the corresponding
field theories. So, the hyperbolic hypergeometric functions are related to supersymmetric partition functions
of three-dimensional field theories. Respectively, complex hypergeometric functions describe
some characteristics of two-dimensional theories and, again, as far as the authors know,
the physical meaning of the additional free variable $\nu=0,\tfrac12$ has not been understood so far.
Therefore we expect some new applications of our results in mathematical physics and hope to
work out some of them in the future.

\appendix

\section{The Mellin-Barnes form of a complex beta integral}

Consider the complex binomial theorem, which is very little discussed in the literature.
It has the following form in the Mellin-Barnes representation
\begin{equation}
\frac{1}{[x+y]^\alpha}=\frac{1}{4\pi {\bf\Gamma}(\alpha|\alpha')}\sum_{N\in\Z}\int_{L} d\nu\,
\frac{{\bf\Gamma}(s|s'){\bf\Gamma}(\alpha-s|\alpha'-s')}{[x]^{\alpha-s}[y]^{s}},
\label{comp_bin}\end{equation}
where $s=\tfrac12(N+\textup{i} \nu)$, $s-s'=N\in\Z$, $\alpha=\tfrac12(m+\textup{i}a)$,
$\alpha-\alpha'=m\in\Z$.
From the integrand factor ${\bf\Gamma}(s|s')=\Gamma(s)/\Gamma(1-s')$ one has formal poles
at $s=-n,\, n\in\Z_{\geq0}$. However, for integer $s$, the variable $s'$ is also an
integer and the same factor has formal zeros if $s'\in\Z_{>0}$.
These zeros should be excluded, since they would cancel the poles, which is reached by imposing
the additional constraint $s'=-n', \, n'\in\Z_{\geq0}$.
As a result, the true pole positions are $\nu= -\textup{i}(s+s')= \textup{i}(n+n')$, which are simple
for fixed $N=s-s'=n' - n$ and go upwards on the complex plane.
Analogously, taking the factor ${\bf\Gamma}(\alpha-s|\alpha'-s')$ and imposing the constraints
$\alpha-s=-n,\, \alpha'-s'=-n'$, one finds the poles in the lower half-plane at
$\nu= \alpha -\textup{i}(n+n')$ for a fixed $N=m+n-n'$.
The integration contour $L$ separates these two sets of poles
and asymptotically lies in any bounded horizontal stripe $|\textrm{Im}(\nu)|<\, const$.
Then the Stirling formula shows that the integrals over $\nu$ converge for Im$(a)>-1$.
If $-1<\mathrm{Im}(a)<0$, then $L$ can be any straight line lying in the stripe
$\mathrm{Im}(a)<\mathrm{Im}(\nu)<0$. In this case, one can shift $s\to s+h$
with $\textrm{Im}(a)< \textrm{Im}(h+h')<0$ and fix $L=\R$, which will be used below.

To prove formula \eqref{comp_bin}, it is necessary to close the integration contour in the upper
half-plane and replace the sum of integrals by the double sum over $n$ and $n'$ of simple
pole residues \cite{Neretin}. Then one can see that
these sums factorize into the product of two $_1F_0$-series converging
 for $|y/x|<1$ and yield the left-hand side expression.
After closing the integration contour in the lower half-plane, the double sum of
residues converges for $|x/y|<1$.
Other values of $x$ and $y$, $x+y\neq 0$, are reached by the analytical continuation.

Relation \eqref{comp_bin} plays an important role in passing from the complex plane integrals
to their Mellin-Barnes form given by the infinite bilateral sums of contour integrals.
For example, replace $z_1-w$ in \eqref{STR} by $x+y$, with $x=z_1-z_3$ and $y=z_3-w$, and
apply formula \eqref{comp_bin} with $\alpha$ replaced by $1-\alpha$. As a result the integral over
$w$ can be computed with the help of the initial identity \eqref{cbeta}.
After simplifications, formula  \eqref{STR} gets reduced to the relation \eqref{comp_bin}
with the replacements $\alpha\to \gamma,\, x\to 1, y\to(z_2-z_2)/(z_1-z_3)$, respectively. This gives
an alternative proof of the star-triangle relation from \eqref{cbeta}
without application of the linear fractional
transformation. Similar mechanism works for passing from the triple complex plane integral describing
$6j$-symbols for the principal series representation of SL$(2,\CC)$ group \cite{Ismag2,DS2017}
to a single infinite sum of univariate Mellin-Barnes type integrals with a bigger number of
complex gamma functions in the integrand than in  \eqref{comp_bin}.

Take relation \eqref{STR} and expand all three $z$-dependent factors on the right-hand side
using \eqref{comp_bin} with $s$ replaced $s+h$ as mentioned above, which yields
$$
\frac{1}{(4\pi)^3}\sum_{N_j\in \Z}\int_{\nu_j\in \R}\frac{\prod_{j=1}^3 d\nu_j\, (-1)^{s_j-s_j'+h-h'}
{\bf \tilde\Gamma}(s_j+h,g_j-h-s_j)}{[z_3]^{g_1-s_1+s_2}[z_2]^{g_3-s_3+s_1}[z_1]^{g_2-s_2+s_3}},
$$
where we denoted $\alpha=g_1, \beta=g_2, \gamma=g_3$ and assumed that
$-1<\textrm{Im}(g_j+g_j')<\textrm{Im}(h+h')<0$.

Multiply now both sides of \eqref{STR} by $[z_1]^{g_2+\rho_1-1}[z_2]^{g_3+\rho_2-1}[z_3]^{g_1+\rho_3-1}$,
where $\rho_j=(k_j+\textup{i}d_j)/2$, $\rho_j-\rho_j'=k_j\in\Z,\, d_j\in\R$, integrate over
all three variables $z_1, z_2, z_3$ using relation \eqref{cbeta}
and apply for $a\in\R$ and complex $z:=re^{i \varphi}$ the standard Fourier transformations relation
$$
\int [z]^{\alpha-1}d(\text{\rm Re}\, z)d(\text{\rm Im}\, z)
=\int_{0}^{2\pi}\int_0^{\infty} e^{\textup{i} \varphi m}r^{ia-2}d\varphi\, rdr
=(2\pi)^2 \delta_{m0}\delta(a)=:(2\pi)^2\delta^{(2)}(\alpha),
$$
where $\delta_{m0}$ and $\delta(a)$ are the discrete and continuous delta functions.

As a result, on the left-hand side one obtains
$$
\pi^2 {\bf\tilde\Gamma}(g_1,g_2+\rho_1,g_3-\rho_1,g_2,g_3+\rho_2,g_1-\rho_2,g_3,
g_1+\rho_3,g_2-\rho_3) I(\underline{\rho}),
$$
where
$$
 I(\underline{\rho})=\int [w]^{\sum_{j=1}^3\rho_j-1}d^2w=(2\pi)^2 \delta^{(2)}(\sum_{j=1}^3\rho_j).
$$

On the right-hand side this procedure yields
$$
\frac{1}{(4\pi)^3}\sum_{N_j\in\Z}\int_{\nu_j\in\R}\prod_{j=1}^3 d\nu_j\, (-1)^{s_j-s_j'+h-h'}\,
{\bf \tilde\Gamma}(s_j+h,g_j-h-s_j)\,\Delta(\underline{s}),
$$
where
$$
\Delta(\underline{s})=(2\pi)^6\delta^{(2)}(\rho_1+s_2-s_3)\delta^{(2)}(\rho_2+s_3-s_1)\delta^{(2)}(\rho_3+s_1-s_2).
$$
Compute sums and integrals over $N_2, N_3, \nu_2, \nu_3$ using the delta functions.
Then, after denoting $s_1=s$, the right-hand side expression takes the form
$$
\pi^3\sum_{N\in\Z}\int_{\R}d\nu\,
\prod_{j=1}^3{\bf\tilde\Gamma}(B_j+s,A_j-s)
(-1)^{s-s'+h-h'+\rho_1-\rho_1'} \delta^{(2)}(\sum_{j=1}^3\rho_j),
$$
where
\begin{eqnarray*} &&
B_1=h,\qquad B_2=\rho_3+h,\qquad B_3= \rho_1+\rho_3+h,
\\ &&
A_1=g_1-h,\qquad A_2=g_2-\rho_3-h,\qquad A_3=g_3-\rho_1-\rho_3-h.
\end{eqnarray*}
Then the left-hand side expression becomes
$$
4\pi^4\prod_{j,k=1}^3{\bf\Gamma}(A_j+B_k|A_j'+B_k')\, \delta^{(2)}(\sum_{j=1}^3\rho_j).
$$
Equating the coefficients in front of the delta functions on both sides under the constraint $\sum_{j=1}^3\rho_j=0$,
one comes \cite{DMV2017} to the required relation \eqref{MB} after setting $A_j=\tfrac12(m_j+\textup{i}a_j),\,
A_j-A_j'=m_j,\, B_j=\tfrac12(n_j+\textup{i}b_j),\, B_j-B_j'=n_j$. The sign factors in the
integrands are removed using the relation ${\bf\Gamma}(x,-n)=(-1)^n{\bf\Gamma}(x,n)$.
It can be checked that the balancing condition \eqref{balancing1} is also satisfied.
The initial assumption that the parameters $\rho_j+\rho_j'$ are purely imaginary
is removed by the analytical continuation.
Finally, we note that it is possible to go from the Mellin-Barnes form of the star-triangle
relation \eqref{MB} back to the integral over complex plane identity \eqref{STR} in the quasiclassical
limit \cite{DM2019}.

\smallskip

{\bf Acknowledgements.}
We would like to thank S. E. Derkachov and A. N. Manashov for valuable discussions and useful comments,
as well as the referees for many constructive remarks and helpful suggestions.
This study has been partially funded by the Russian Science Foundation (grant 24-21-00466).

\end{document}